\documentclass{article}%
\usepackage{epsfig}
\usepackage{enumerate,graphicx}
\usepackage{amsfonts}
\usepackage{graphicx}
\usepackage{amsmath,amsthm,amssymb}
\usepackage{amssymb}%
\DeclareGraphicsExtensions{.pdf, .jpg, .png , .bmp, .ps}
\usepackage{tikz}
\providecommand{\U}[1]{\protect\rule{.1in}{.1in}}
\newtheorem{thm}{Theorem}[section]

\newtheorem{prop}{Proposition}[section]

\theoremstyle{remark}
\newtheorem{rem}[thm]{\bf Remark}

\def\1{{{\mbox{${\rm{1\negthinspace\negthinspace I}}$}}}}

\newcommand{\eref}[1]{(\ref{#1})}

\newcommand\beq{\begin{equation}}
\newcommand\eeq{\end{equation}}

\begin{document}

\title{McDiarmid's martingale for a class of iterated random functions}

\author{J\'er\^ome Dedecker\footnote{Universit\'e Paris Descartes,
Sorbonne Paris Cit\'e,
Laboratoire MAP5
and CNRS UMR 8145, 75016 Paris, France.}\footnotemark[1]\ \ \ and \  Xiequan Fan\footnote{Regularity Team, Inria and MAS Laboratory, Ecole Centrale Paris - Grande Voie des Vignes, 92295 Ch\^{a}tenay-Malabry, France.}\footnotemark[2]}

\date{}

\maketitle

\abstract{We consider an ${\mathcal X}$-valued Markov chain $X_1, X_2, \ldots, X_n$  belonging to a class of iterated random functions, which is
``one-step contracting" with respect to some distance $d$ on ${\mathcal X}$. If $f$ is any  separately Lipschitz
function with respect to $d$, we use a well known decomposition of
$S_n=f(X_1, \ldots, X_n) -{\mathbb E}[f(X_1, \ldots , X_n)]$ into a sum of martingale
differences $d_k$
with respect to the  natural filtration ${\mathcal F}_k$. We show that each difference $d_k$
is
bounded by a random variable $\eta_k$ independent of ${\mathcal F}_{k-1}$.
Using this very strong property, we obtain a large variety of deviation inequalities
for $S_n$, which are governed by the distribution of the $\eta_k$'s.
Finally, we give an application of these inequalities to the Wasserstein distance
between the empirical measure and the invariant distribution of the chain.
 }

\medskip


\medskip


\section{A class of iterated random functions}

Let $(\Omega, {\mathcal A}, {\mathbb P})$ be a probability space.
Let $({\mathcal X}, d)$ and  $({\mathcal Y}, \delta)$ be two complete separable metric spaces.
Let $(\varepsilon_i)_{i \geq 1}$ be a sequence of independent and identically distributed (iid) ${\mathcal Y}$-valued
random variables. Let $X_1$ be a ${\mathcal X}$-valued random variable independent of $(\varepsilon_i)_{i \geq 2}$. We consider the Markov chain $(X_i)_{i \geq 1}$
such that
\begin{equation}\label{Mchain}
X_n=F(X_{n-1}, \varepsilon_n), \quad \text{for $n\geq 2$},
\end{equation}
where $F: {\mathcal X} \times {\mathcal Y} \rightarrow {\mathcal X}$
is such that
\begin{equation}\label{contract}
{\mathbb E}\big [ d\big(F(x, \varepsilon_1), F(x', \varepsilon_1)\big) \big] \leq
\rho d(x, x')
\end{equation}
for some $\rho \in [0,1)$, and
\begin{equation}\label{c2}
d(F(x,y), F(x,y')) \leq C \delta(y,y')
\end{equation}
for some  $C>0$.

This class of Markov chains, that we call ``one-step contracting", is very restrictive, but still contains a lot
of pertinent examples. Among them, in the case where ${\mathcal X}$ is a separable Banach space with norm $|\cdot|$, let us cite the functional auto-regressive model
$$
X_n=f(X_{n-1}) + g(\xi_n)\, ,
$$
where $f : {\mathcal X}  \rightarrow {\mathcal X}$ and $g: {\mathcal Y}  \rightarrow {\mathcal X}$ are such that
$$
  |f(x)-f(x')|\leq \rho |x-x'|  \quad \text{and} \quad |g(y)-g(y')| \leq C \delta(y,y') \, .
$$
We refer to the paper by Diaconis and Freedman \cite{DF99} for many other interesting examples.
Note also that this class of Markov chains contains the iid sequence $X_i=\varepsilon_i$, by taking
${\mathcal Y}={\mathcal X}$  and $F(x,y)=y$ (note that $\rho=0$ in that case).

This class possesses the property of exponential forgetting of the starting point: If $X_n^x$ is the chain
starting from $X_1=x$, then one has
$$
  {\mathbb E}\big[d(X_n^x, X_n^{x'})\big] \leq \rho^n d(x,x')\, .
$$
Hence is has an unique stationary distribution $\mu$ (see for instance Theorem 1 in Diaconis and Freedman \cite{DF99}), meaning
that if $X_1$ is distributed as $\mu$, then the chain $(X_i)_{i \geq 1}$ is strictly stationary. Moreover, one can easily
prove that, if $(X_i)_{i \geq 1}$ is strictly stationary, then, for any $(x_0, y_0) \in {\mathcal X} \times {\mathcal Y}$,
and any positive measurable function $H$,
\begin{equation}\label{trivial}
{\mathbb E}[H(d(X_n,x_0))] \leq {\mathbb E} \Big[H\Big(\sum_{i=0}^\infty \rho^i \big(d(F(x_0,y_0), x_0) + C \delta(\varepsilon_{i+1}, y_0)\big)\Big)\Big]\, .
\end{equation}

Although the one-step contraction is a very restrictive condition, this class of iterated
random functions contains a lot
of non Harris-recurrent Markov chains. For instance, if ${\mathcal X}={\mathcal Y}=[0,1]$ the chain
$$
X_{n}=\frac{1}{2}(X_{n-1}+ \varepsilon_n)
$$
with $X_1$ uniformly distributed over $[0,1]$, and $\varepsilon_i \sim {\mathcal B}(1/2)$ is strictly
stationary, but it is not mixing in the sense of Rosenblatt \cite{Ro}.

The class of iteretad random function satsifying \eref{contract} has been studied in
Section 3.1 of
Djellout {\it et al.} \cite{DGW} (as a particular case of a general class of Markov chains
which are contracting with respect to Wasserstein distances, see their Condition C1).
Combining  McDiarmid method and a result by Bobkov and G\"otze \cite{BoG}, Djellout {\it et al.}
\cite{DGW}
proved in their Proposition 3.1   a subgaussian bound for
separately Lipschitz
functionals of the chain provided
\begin{equation}\label{DGWC}
\sup_{x \in {\mathcal X}} {\mathbb E}
\Big [ \exp \Big( a \big(d(F(x,\varepsilon_1), F(x, \varepsilon_2))\big)^2 \Big) \Big]
< \infty \, ,
\end{equation}
for some $a>0$.
Because of the supremum in $x$, this condition is quite delicate to check. However,
if \eref{c2} holds, it is implied by the simple condition
$$
{\mathbb E}
\Big [ \exp \Big( a \big(C \delta(\varepsilon_1,  \varepsilon_2)\big)^2 \Big) \Big] < \infty \, .
$$
As we shall see in Section \ref{secMcD}, this is due to the fact that the martingale
differences from McDiarmid's decomposition are bounded by a random variable
$\eta_k$ independent of
${\mathcal F}_{k-1}= \sigma (X_1, \ldots, X_{k-1})$. From this simple remark,
we can obtain many deviation inequalities for separately Lipschitz functionals of the chain
by applying known inequalities for martingales.

 A more restrictive class of iterated random function, satisfying \eref{c2} and the one-step contraction
 $$
  d\big(F(x, y), F(x', y)\big)  \leq
\rho d(x, x')\, ,
 $$
 has been studied by Delyon {\it et al.} \cite{DJL06} when
${\mathcal X}= {\mathbb R}^\ell$ and ${\mathcal Y}={\mathbb R}^k$. These authors have proved a moderate
deviation principle for additive and Lipschitz functionals of the chain, under a condition on the
Laplace transform of the euclidean norm of $\varepsilon_i$.



\section{McDiarmid's martingale}\label{secMcD}
\setcounter{equation}{0}
\subsection{Separately Lipschitz functions of $X_1, \ldots , X_n$.}
Let $f: {\mathcal X}^n \mapsto {\mathbb R}$ be separately Lipschitz, such that
\begin{equation} \label{codiMD}
|f(x_1, x_2, \ldots, x_n)-f(x'_1, x'_2, \ldots, x'_n)| \leq d(x_1,x'_1)+\cdots +  d(x_n, x'_n) \, .
\end{equation}
Let then
\begin{equation}\label{Sn}
S_n:=f(X_1, \ldots , X_n) -{\mathbb E}[f(X_1, \ldots , X_n)]\, .
\end{equation}
We also introduce the natural filtration of the chain, that is  ${\mathcal F}_0=\{\emptyset, \Omega \}$
and for $k \in {\mathbb{N}}^{*}$,
${\mathcal F}_k= \sigma(X_1, X_2,  \ldots, X_k)$.
Define then
\begin{equation}\label{gk}
g_k(X_1, \ldots , X_k)= {\mathbb E}[f(X_1, \ldots, X_n)|{\mathcal F}_k]\, ,
\end{equation}
and
\begin{equation}\label{dk}
d_k=g_k(X_1, \ldots, X_k)-g_{k-1}(X_1, \ldots, X_{k-1})\, .
\end{equation}
For $k \in [1, n-1]$, let
$$
S_k:=d_1+d_2+\cdots + d_k \, ,
$$
and note that, by definition of the $d_k$'s, the functional $S_n$ introduced in
\eref{Sn} satisfies
$$
S_n =d_1+d_2+\cdots + d_n \, .
$$
Hence $S_k$ is a martingale adapted to the filtration ${\mathcal F}_k$. This representation was introduced by McDiarmid \cite{M} in
the iid case, when $X_i=\varepsilon_i$ (see also Yurinskii \cite{Yu} in a different context).

The following Proposition collects some interesting properties of
the functions $g_k$ and of the martingale differences $d_k$.

\begin{prop}\label{McD}
For $k \in {\mathbb N}$
and $\rho$ in $[0, 1)$, let  $K_k(\rho)=(1-\rho^{k+1})/(1-\rho)= 1 + \rho + \cdots + \rho^{k}$.
Let $(X_i)_{i \geq 1}$ be a Markov chain satisfying \eref{Mchain} for some
function $F$ satisfying \eref{contract}. Let $g_k$ and $d_k$ be defined
by
\eref{gk} and \eref{dk} respectively.
\begin{enumerate}
\item
The function $g_k$ is separately Lipschitz and such that
$$
|g_k(x_1, x_2, \ldots, x_k)-g_k(x'_1, x'_2, \ldots, x'_k)| \leq d(x_1,x'_1)+\cdots + d(x_{k-1}, x'_{k-1})+ K_{n-k}(\rho) d(x_k, x'_k) \, .
$$
\item Let $P_{X_1}$ be the distribution of $X_1$ and $P_\varepsilon$ be the common
distribution of the $\varepsilon_k$'s. Let $G_{X_1}$ and $H_{\varepsilon}$ be the two functions defined by
$$
G_{X_1}(x)=\int d(x, x') P_{X_1}(dx') \quad \text{and}
\quad
H_{\varepsilon}(x,y)= \int d(F(x,y),F(x,y'))P_{\varepsilon}(dy') \, .
$$
Then, the martingale difference $d_k$ is such that
$$
   |d_1| \leq K_{n-1}(\rho) G_{X_1}(X_1) \quad
   \text{and for $k \in [2,n]$,} \quad
   |d_k|\leq K_{n-k}(\rho)  H_\varepsilon(X_{k-1}, \varepsilon_k)\, .
$$
\item
Assume moreover that $F$ satisfies \eref{c2}, and
Let $G_{\varepsilon}$ be the  function defined by
$$
G_{\varepsilon}(y)= \int C\delta(y,y')P_{\varepsilon}(dy') \, .
$$
Then $H_\varepsilon(x,y) \leq G_\varepsilon(y)$,
and consequently, for $k \in [2,n]$,
$$
 |d_k|\leq K_{n-k}(\rho)  G_\varepsilon(\varepsilon_k)\, .
$$
\end{enumerate}
\end{prop}

\begin{rem}
Let us comment on the point 3 of Proposition \ref{McD}.
The fact that the martingale difference $d_k$ is bounded by the random variable
$K_{n-k}(\rho)G_\varepsilon(\varepsilon_k)$
which is {\it independent  of }${\mathcal F}_{k-1}$ is crucial. It explains
why we shall obtain deviations inequalities for $S_n$ under some conditions on the
distribution of $G_\varepsilon(\varepsilon_k)$ (typically conditions on the Laplace
transform, or moment conditions).
\end{rem}

\medskip

\noindent {\emph {Proof}.} The first point will be proved by recurrence in the backward sense.
The result is obvious for $k=n$, since $g_n=f$. Assume that it is true at step $k$, and let us prove it
at step $k-1$. By definition
$$
g_{k-1}(X_1, \ldots, X_{k-1})={\mathbb E}[g_k(X_1, \ldots, X_k)|{\mathcal F}_{k-1}]= \int g_k(X_k, \ldots, X_{k-1}, F(X_{k-1},y)) P_{\varepsilon}(dy)\, .
$$
It follows that
\begin{multline}\label{triv1}
|g_{k-1}(x_1, x_2, \ldots, x_{k-1})-g_{k-1}(x'_1, x'_2, \ldots, x'_{k-1})|\\ \leq \int |g_k(x_1, x_2, \ldots, F(x_{k-1},y))-g_k(x'_1, x'_2, \ldots, F(x'_{k-1},y))|
P_{\varepsilon}(dy)\, .
\end{multline}
Now, by assumption and condition \eref{contract},
\begin{multline}\label{triv2}
\int |g_k(x_1, x_2, \ldots, F(x_{k-1},y))-g_k(x'_1, x'_2, \ldots, F(x'_{k-1},y))|
P_{\varepsilon}(dy)\\
\leq d(x_1,x'_1)+\cdots + d(x_{k-1}, x'_{k-1})+ K_{n-k}(\rho) \int d(F(x_{k-1},y), F(x'_{k-1},y))P_{\varepsilon}(dy)
\\
\leq d(x_1,x'_1)+\cdots  + (1+\rho K_{n-k}(\rho)) d(x_{k-1}, x'_{k-1}) \\
\leq d(x_1,x'_1)+\cdots  +  K_{n-k+1}(\rho) d(x_{k-1}, x'_{k-1}) \, .
\end{multline}
The point 1 follows from \eref{triv1} and \eref{triv2}.

Let us prove the point 2. First note that
$$
|d_1|=\Big|g_1(X_1)-\int g_1(x)P_{X_1}(dx)\Big|\leq
 K_{n-1}(\rho)\int d(X_1,x)P_{X_1}(dx)=K_{n-1}(\rho) G_{X_1}(X_1) \, .
$$
In the same way, for $k \geq 2$,
\begin{align*}
 |d_k| & = \big|g_k(X_1, \cdots, X_k)-{\mathbb E}[g_k(X_1, \cdots, X_k)|{\mathcal F}_{k-1}]\big|\\
 &\leq  \int \big |g_k(X_1, \cdots, F(X_{k-1}, \varepsilon_k))-g_k(X_1, \cdots, F(X_{k-1}, y))\big| P_{\varepsilon}(dy)\\
 &\leq K_{n-k}(\rho)\int  d(F(X_{k-1},\varepsilon_k),
 F(X_{k-1}, y))P_{\varepsilon}(dy)= K_{n-k}(\rho)
  H_{\varepsilon}(X_{k-1},\varepsilon_k) \, .
\end{align*}

The point 3 is clear, since if \eref{c2} is true, then
$$
H_{\varepsilon}(x,y)= \int d(F(x,y),F(x,y'))P_{\varepsilon}(dy')
\leq \int C\delta(y,y')P_{\varepsilon}(dy')= G_\varepsilon(y) \, .
$$
The proof of the proposition is now complete.
\hfill\qed

\subsection{An important remark}
For any $\alpha \in (0,1)$ define the distances $d_{\alpha}$ and $\delta_\alpha$ on
${\mathcal X}$ and ${\mathcal Y}$ respectively by
$$
  d_\alpha(x,x')= (d(x,x'))^\alpha \quad \text{and} \quad
  \delta_\alpha(y,y') = (\delta(y,y'))^\alpha \, .
$$

If $F$ is one-step contacting with respect to a natural distance $d$ (meaning
that it satisfies the inequalities \eref{contract}  and \eref{c2}
with $\rho \in [0,1)$ and $C>0$ respectively), then
for any $\alpha \in (0,1)$,
\begin{equation}\label{contractbis}
{\mathbb E}\big[ \big(d_\alpha(F(x,\varepsilon_1), F(x',\varepsilon_1)\big) \big]
\leq \rho^\alpha d_\alpha(x, x') \, ,
\end{equation}
and
\begin{equation}
d_\alpha \big (F(x,y), F(x, y')\big) \leq C^\alpha
\delta_\alpha(y,y') \, .
\end{equation}
Hence $F$ is also one-step contracting for the distance $d_\alpha$, with the new
constants $\rho^\alpha \in [0,1)$ and $C^\alpha>0$.
Consequently, Proposition \ref{McD} applies to the martingale
$$
S_n=f(X_1, \ldots , X_n) -{\mathbb E}[f(X_1, \ldots , X_n)]\, ,
$$
where $f$ is separately Lipshitz with respect to $d_\alpha$. The dominating random variables $G_{X_1, \alpha}(X_1)$ and $G_{\varepsilon, \alpha}(\varepsilon_k)$ are then defined by
$$
G_{X_1, \alpha}(x)=\int d_\alpha(x, x') P_{X_1}(dx') \quad \text{and}
\quad\,
G_{\varepsilon, \alpha}(y)= \int C^\alpha \delta_\alpha(y,y')P_{\varepsilon}(dy')  \, .
$$
Hence, all the results of the following section apply to the functional $S_n$,
provided the corresponding conditions on the dominating random  variables
$G_{X_1, \alpha}(X_1)$ and $G_{\varepsilon, \alpha}(\varepsilon_k)$ are satisfied.

For instance, if ${\mathcal X}={\mathbb R}^\ell$ and $d(x,y)=\|x-y\|$ is the
euclidean distance on ${\mathbb R}^\ell$, then one can consider the class
of
separately H\"older functions $f$ such that
\begin{equation*} \label{codiMD2}
|f(x_1, x_2, \ldots, x_n)-f(x'_1, x'_2, \ldots, x'_n)| \leq \|x_1 - x'_1\|^\alpha+
\cdots +  \|x_n - x'_n\|^\alpha \, .
\end{equation*}

\section{Deviation inequalities for the functional $S_n$.} \label{deviationiq}

\setcounter{equation}{0}
Let $(X_i)_{i \geq 1}$ be a Markov chain satisfying \eref{Mchain} for some
function $F$ satisfying \eref{contract} and \eref{c2}.
In this section, we apply inequalities for martingales to bound up
the deviation of the functional $S_n$ defined by \eref{Sn}. Some of these inequalities
are direct applications of known inequalities, some deserve a short proof and some
other are new.

Note that deviation inequalities for Lipschitz functions of dependent sequences
have been proved for instance by Rio \cite{R00}, Collet {\it et al.} \cite{CMS},
Djellout {\it et al.} \cite{DGW}, Kontorovich and Ramanan \cite{KR08},  and Chazottes and Gou\"ezel \cite{CG12} among others.
Except for Djellout {\it et al.} \cite{DGW} (who also consider more general Markov chains),
the examples  studied by these authors are different from the class described in the present paper. For instance, the Markov chains associated to the
maps studied by Chazottes and Gou\"ezel \cite{CG12} do not in general satisfy the one step
contraction property.

The interest of the one step contraction is that, thanks to Proposition \ref{McD}, we shall obtain very precise inequalities,
with precise constants depending on the distribution of the dominating random
variables $G_{X_1}(X_1)$ and $G_\varepsilon(\varepsilon_k)$.

Let us note that, in the iid case, when $X_i=\varepsilon_i$, the additive functional
$$
f(x_1, x_2, \ldots, x_n)= \sum_{k=1}^n G_{\varepsilon}(x_i)
$$
is of course separately Lipshitz and satisfies \eref{codiMD}.
Hence, the inequalities of the following section apply to this
simple functional, under the usual moment or Laplace conditions
on the (non centered) variables $G_{\varepsilon}(\varepsilon_i)$.
This shows that, in the iid case, these inequalities cannot be much improved
without additional assumptions on the functional $f$.

Les us now consider the case where we only assume that $F$ satisfies \eref{contract}.
Then all the inequalities of this section will be true provided the appropriate conditions
of the type
${\mathbb E}[f(G_\varepsilon(\varepsilon))]\leq C$ for some positive measurable function
$f$ are replaced by
\begin{equation}\label{weaker}
\sup_{k \in [2,n]}
 \Big \| {\mathbb E} \Big [ f\big ( H_\varepsilon(X_{k-1}, \varepsilon_k) \big)
 \Big | X_{k-1}\Big ]\Big \|_\infty \leq C \, .
\end{equation}
Note that the latter condition is true provided
$$
\sup_{x \in {\mathcal X}}
  {\mathbb E} \Big [ f\big ( H_\varepsilon(x, \varepsilon_1) \big)
 \Big ] \leq C \, ,
$$
which is of the same type as condition \eref{DGWC}
for the subgaussian bound (with $f(x)=\exp(ax^2)$ in that particular case).
Recall that condition \eref{DGWC} is due to  Djellout {\it et al.} \cite{DGW} (see their
Proposition 3.1).

For the weak and strong moment bounds on $S_n$, we shall see in Subsections \ref{VBEB}, \ref{MZB} and \ref{BRB} that
 condition \eref{weaker} can be replaced by an appropriate moment condition on
$H_\varepsilon(X_{k-1}, \varepsilon_k)$.

To conclude the introduction of this section, let us note that the deviations
inequalities of Subsections \ref{Bersec} -- \ref{Fuksec} are given for ${\mathbb P}\big(\pm S_n >x\big)$,
but we shall only prove them for $S_n$. The proofs of the
deviation inequalities for $-S_n$ are exactly the same, the upper bounds of  points 2 and 3 of Proposition \ref{McD} being valid for  $d_k$ and $-d_k$.

In all this section,  $G_\varepsilon(\varepsilon)$ denotes
a random variable
distributed as $G_\varepsilon(\varepsilon_k)$.


\subsection{Bernstein type bound}\label{Bersec}
Under the conditional Bernstein condition, van de Geer \cite{V95} and De La Pe\~{n}a \cite{D99} have obtained some tight Bernstein type inequalities for  martingales.
Applying Proposition \ref{McD}, we obtain the following proposition.
\begin{prop}
Assume that there exist some constants  $M>0, V_1\geq 0$ and $V_2\geq 0$ such that, for any integer $k\geq 2$,
\begin{equation}\label{Bernsteinmoment}
{\mathbb E} \Big[  \Big(  G_{X_1}(X_1)\Big)^k\Big] \leq \frac {k!}{2} V_1 M^{k-2}
\ \text{and} \quad
 {\mathbb E} \Big[  \Big(  G_\varepsilon(\varepsilon)\Big)^k\Big] \leq
 \frac {k!}{2} V_2 M^{k-2} \, .
\end{equation}
Let $$V=  V_1\Big(K_{n-1}(\rho)\Big)^2 +  V_2 \sum_{k=2}^n \Big(K_{n-k}(\rho)\Big)^2  \ \ \ \ \ \textrm{and}\ \ \ \ \ \delta=MK_{n-1}(\rho).$$
Then, for any $t \in [0, \delta^{-1})$,
\begin{equation}\label{maindfs}
  \mathbb{E}\,[e^ {\pm  tS_n} ]\leq \exp \left (\frac{t^2 V }{2 (1- t\,\delta )} \right )\, .
\end{equation}
Consequently, for any $x> 0$,
\begin{eqnarray}
  {\mathbb P}\big(\pm  S_n \geq  x\big)
  &\leq&  \exp \left( \frac{x^2}{V(1+\sqrt{1+2x \delta/V})+x \delta }\right)\, \label{Berie1} \\
 &\leq&  \exp \left( \frac{x^2}{2 \,(V +x  \delta ) }\right)\, .\label{Berie2}
\end{eqnarray}
\end{prop}
\begin{rem}
Let us comment on condition \eref{Bernsteinmoment}.
\begin{enumerate}
\item In the iid case, when $X_i=\varepsilon_i$, condition \eref{Bernsteinmoment}
is the Bernstein condition
$$
{\mathbb E} \Big[  \Big(   G_\varepsilon(\varepsilon)\Big)^k\Big] \leq
 \frac {k!}{2} V M^{k-2}.
$$
In that case the inequalities \eref{Berie1} and \eref{Berie2} hold with $\rho=0$.
\item
Since
$G_\varepsilon(\varepsilon) \leq C\delta(\varepsilon, y_0) +
C{\mathbb E}[\delta(\varepsilon, y_0)]$, it follows that
$$
{\mathbb E} \Big[ \Big( G_\varepsilon(\varepsilon)\Big)^k\Big] \leq 2^k {\mathbb E} \Big[ \Big( C\delta(\varepsilon, y_0)\Big)^k\Big]\, .
$$
Hence,
the condition
\begin{equation}\label{Berbis}
{\mathbb E} \Big[ \Big( C\delta(\varepsilon, y_0)\Big)^k\Big] \leq
\frac {k!}{2} A(y_0) B(y_0)^{k-2}
\end{equation}
implies the second condition in \eref{Bernsteinmoment}
with $V_2=4A(y_0)$ and $M=2B(y_0)$.
In the same way, the condition
\begin{equation}\label{Berter}
{\mathbb E} \Big[ \Big( d(X_1, x_0)\Big)^k\Big] \leq
\frac {k!}{2} C(x_0) D(x_0)^{k-2}
\end{equation}
implies the first condition in \eref{Bernsteinmoment}
with $V_1=4C(x_0)$ and $M=2D(x_0)$.
\item  Consider the chain with non random starting point $X_1=x$. Then
$G_{X_1}(X_1) = 0$, and the first condition in \eref{Bernsteinmoment}
holds with $V_1=0$.
\item Let us consider now the case where $X_1$ is distributed according to the
invariant probability measure $\mu$.
We shall see that in that case \eref{Berter}
follows from \eref{Berbis}.
To avoid to many computations, assume that one can find
$(x_0, y_0)$ such that
$d(F(x_0,y_0), x_0)=0$, which is true in many cases.
If \eref{Berbis} holds, it follows from \eref{trivial} applied to $H(x)=x^k$ that
\eref{Berter} holds with $C(x_0)=(1-\rho)^{-2}A(y_0)$ and $D(x_0)=(1-\rho)^{-1}B(y_0)$.
According to the point 2 of this remark, condition \eref{Bernsteinmoment} is satisfied by taking
$M=2(1-\rho)^{-1}B(y_0)$, $V_2=4A(y_0)$ and $V_1=4(1-\rho)^{-2}A(y_0)$.
\end{enumerate}
\end{rem}
\noindent\emph{Proof.} From Proposition \ref{McD} and condition (\ref{Bernsteinmoment}), it is  easy to see that, for any $t \in [0, \delta^{-1})$,
\begin{eqnarray}\label{fin35}
\mathbb{E}\,[e^{t d_1 } ] &\leq& 1+ \sum_{i=2}^{\infty} \frac{t^i}{i!}\, \mathbb{E}\,[   (d_1)^i ] \nonumber\\
&\leq& 1+ \sum_{i=2}^{\infty} \frac{t^i}{i!}\, \mathbb{E}\,[ | d_1|^i ] \nonumber\\
&\leq& 1+ \sum_{i=2}^{\infty} \frac{t^i}{i!}\, \Big(K_{n-1}(\rho)\Big)^i \mathbb{E}\,\Big[ \Big(G_{X_1}(X_1)\Big)^i \Big] \nonumber\\
&\leq& 1+ \sum_{i=2}^{\infty} \frac{t^i}{i!}\, \Big(K_{n-1}(\rho)\Big)^i \frac {i!}{2} V_1 M^{i-2}
= 1+ \frac{t^2V_1 \Big(K_{n-1}(\rho)\Big)^2  }{2 (1  -t\, \delta )}\, .
\end{eqnarray}
Similarly,   for any  $k \in [2, n],$
\begin{eqnarray}\label{fin36}
\mathbb{E}\,[e^{t d_k  }|\mathcal{F}_{k-1}] \leq 1+ \frac{t^2V_2\Big(K_{n-k}(\rho)\Big)^2}{2 (1  -t\, \delta )} \, .
\end{eqnarray}
Using the inequality $1+t \leq e^t,$ we find that, for any $t \in [0, \delta^{-1})$,
\begin{eqnarray}
\mathbb{E}\,[e^{t d_1 } ]
&\leq&  \exp \left (\frac{t^2V_1\Big(K_{n-1}(\rho)\Big)^2}{2 (1  -t\,\delta )}  \right )
\end{eqnarray}
and
\begin{eqnarray}
\mathbb{E}\,[e^{t d_k }|\mathcal{F}_{k-1}] &\leq& \exp \left( \frac{t^2V_2\Big(K_{n-k}(\rho)\Big)^2}{2 (1  -t\,\delta)}  \right).
\end{eqnarray}
By the tower property of conditional expectation, it follows that, for any $t \in [0, \delta^{-1})$,
\begin{eqnarray*}
\mathbb{E}\,\big[e^{ tS_n} \big] &=&  \mathbb{E}\,\big[ \mathbb{E}\, [e^{ tS_n} |\mathcal{F}_{n-1}  ] \big]\\
&=&  \mathbb{E}\,\big[ e^ { tS_{ n-1}} \mathbb{E}\,  [e^ { td_n}  |\mathcal{F}_{n-1}  ] \big] \\
&\leq & \mathbb{E}\,\big[ e^ { tS_{ n-1}} \big] \exp \left( \frac{t^2V_2}{2 (1- t\,\delta)}  \right)\\
&\leq &  \exp \left(\frac{t^2 V}{2 (1- t\,\delta)} \right),
\end{eqnarray*}
which gives inequality (\ref{maindfs}).
Using the exponential Markov inequality, we deduce that,
for any $x\geq 0$ and $t \in [0, \delta^{-1})$,
\begin{eqnarray}
  \mathbb{P}\left( S_{n} \geq x \right)
 &\leq&  \mathbb{E}\, \big[e^{t\,(S_n -x) } \big ] \nonumber\\
 &\leq&   \exp \left(-t\,x  +  \frac{t^2 V}{2 (1- t\,\delta)}    \right)\, . \label{fines}
\end{eqnarray}
The minimum is reached at $$t=t(x):= \frac{2x/V}{ 2x\delta/V +1 + \sqrt{1+2x\delta/V} }\, .$$ Substituting $t=t(x)$ in (\ref{fines}),  we obtain the desired inequalities
\begin{eqnarray}
  \mathbb{P}\left( S_{n} \geq x \right)
 &\leq&  \exp \left( \frac{x^2}{V(1+\sqrt{1+2x \delta/V})+x  \delta }\right)\, \nonumber\\
 &\leq&  \exp \left( \frac{x^2}{2 (V +x  \delta ) }\right)\, ,\nonumber
\end{eqnarray}
where the last line follows from the inequality $\sqrt{1+2x\,\delta/V}  \leq 1+x\, \delta/V$. \hfill\qed


\subsection{Cram\'{e}r type bound}
If the Laplace transform of the dominating random variables $G_{X_1}(X_1)$ and $G_{\varepsilon}(\varepsilon_k)$ satisfy the Cram\'{e}r condition, we obtain the following proposition
similar to that of Liu and Watbled \cite{LW09}  under the conditional Cram\'{e}r condition. For the optimal convergence speed of martingales under the Cram\'{e}r condition, we refer to Lesigne of Voln\'{y} \cite{LV01} and  Fan {\it et al.} \cite{Fx1}.
\begin{prop}\label{cram}
Assume that there exist some constants $a>0,$ $K_1\geq 1$ and $K_2\geq 1$ such that
\begin{equation}\label{laplace}
{\mathbb E} \Big[ \exp \Big( a G_{X_1}(X_1)\Big)\Big] \leq K_1 \quad \text{and} \quad
 {\mathbb E} \Big[ \exp \Big( a C G_\varepsilon(\varepsilon)\Big)\Big] \leq K_2
 \, .
\end{equation}
Let $$K=\frac{2}{e^{2}} \left( K_1+ K_2\sum_{i=2}^{n} \Big(\frac {K_{n-i}(\rho)}{K_{n-1}(\rho)}\Big)^2 \right)  \ \ \ \ \ and \ \ \ \ \ \delta=\frac{a}{K_{n-1}(\rho)} .$$
Then, for any $t \in [0, \delta )$,
$$
  \mathbb{E}\,[e^{\pm  tS_n}]\leq   \exp \left( \frac{t^2 K  \delta^{-2} }{1-t \delta^{-1}}   \right) \, .
$$
Consequently, for any $x> 0$,
\begin{eqnarray}
  {\mathbb P}\big(\pm  S_n \geq  x\big) &\leq&  \exp \left( \frac{(x\delta)^2}{2K (1+\sqrt{1+ x \delta / K })+x \delta  }\right)\,   \label{Berie3}  \\
 &\leq&  \exp \left( \frac{(x\delta)^2}{ 4K+ 2 x  \delta  }\right)\,\, .\label{Berie4}
\end{eqnarray}
\end{prop}
\begin{rem}
Let us comment on condition \eref{laplace}.
\begin{enumerate}
\item In the iid case, when $X_i=\varepsilon_i$, the condition \eref{laplace}
writes simply
$$
{\mathbb E} \Big[ \exp \Big( aG_{\varepsilon}(\varepsilon)\Big)\Big] \leq K \, .
$$
In that case the inequalities \eref{Berie3} and \eref{Berie4} hold with $\rho=0$.
\item
Since
$G_\varepsilon(\varepsilon) \leq C\delta(\varepsilon, y_0) +
C{\mathbb E}(\delta(\varepsilon, y_0))$
the condition
\begin{equation}\label{bis}
{\mathbb E} \Big[ \exp \Big( a C\delta(\varepsilon, y_0)\Big)\Big] \leq A(y_0)
\end{equation}
implies the second condition in \eref{laplace}
with $K_2=A(y_0)
\exp \big( a C {\mathbb E}[\delta(\xi, y_0)]\big)\leq A(y_0)^2$.
In the same way, the condition
\begin{equation}\label{ter}
{\mathbb E} \Big[ \exp \Big( a d(X_1, x_0)\Big)\Big] \leq B(x_0)
\end{equation}
implies the first condition in \eref{laplace}  with $K_1=B(x_0)
\exp \big(a{\mathbb E}[d(X_1, x_0)]\big)\leq B(x_0)^2$.
\item Consider the chain with non random starting point $X_1=x$. Then
$G_{X_1}(X_1) = 0$, and the first condition in \eref{laplace}
holds with $K_1=1$.
\item Let us consider now the case where $X_1$ is distributed according to the
invariant probability measure $\mu$. We shall see that in that case \eref{ter}
follows from \eref{bis}. Indeed, if \eref{bis} holds, it follows from \eref{trivial} applied
to $H(x)=\exp(a x )$ that
$$
{\mathbb E} \Big[ \exp \Big( a d(X_1, x_0)\Big)\Big] \leq
\exp \Big( \frac{a}{1-\rho}d(F(x_0,y_0), x_0) \Big) \prod_{i=0}^\infty (A(y_0))^{\rho^i}\, .
$$
Hence
$$
{\mathbb E} \Big[ \exp \Big( a d(X_1, x_0)\Big)\Big] \leq
\exp \Big( \frac{a}{1-\rho}d(F(x_0,y_0), x_0) \Big) (A(y_0))^{1/(1-\rho)}
$$
and  \eref{ter} is true with
$$B(x_0)=(A(y_0))^{1/(1-\rho)}\exp \Big( \frac{a}{1-\rho}d(F(x_0,y_0), x_0) \Big)\, . $$
According to the point 2 of this remark, condition \eref{laplace} is satisfied
by taking
$K_2= (A(y_0))^2$ and $K_1= (B(x_0))^2$.
In particular, if \eref{bis} holds, and if we can find $(x_0, y_0)$ such that
$d(F(x_0,y_0), x_0)=0$, then one can take
$K_1=(A(y_0))^{2/(1-\rho)}$.
\end{enumerate}
\end{rem}

\noindent\emph{Proof.} Let $\delta=a/K_{n-1}(\rho).$  Since $\mathbb{E}\,[  d_1 ] =0$, it is  easy to see that, for any $t \in   [0, \delta )$,
\begin{eqnarray}\label{finsa3f}
\mathbb{E}\,[e^{t d_1 } ] &=& 1+ \sum_{i=2}^{\infty} \frac{t^i}{i!}\, \mathbb{E}\,[   (d_1)^i ] \nonumber\\
&\leq& 1+ \sum_{i=2}^{\infty} \Big(\frac{t}{\delta} \Big)^i  \, \mathbb{E}\,\Big[\frac{1}{i!} | \delta  d_1|^i \Big].
\end{eqnarray}
Here, let us note that, for $t\geq 0$,
\begin{align}\label{constantsecond}
\frac{t^i}{i !} e^{-t} &\leq \frac{i^i}{i !} e^{-i} \nonumber \\
&\leq 2 e^{-2}, \quad \text{for $i\geq 2$,}
\end{align}
where the last line follows from the fact that $i^i e^{-i}/i!$ is decreasing in $i$. Note that the equality in \eref{constantsecond}
 is reached at $t=i=2$. Using \eref{constantsecond}, Proposition \ref{McD} and condition (\ref{laplace}), we have
\begin{eqnarray}\label{fisa3f}
\mathbb{E}\,\Big[\frac{1}{i!} |\delta d_1|^i \Big]   &\leq& \ 2e^{-2}  \mathbb{E}\, [e^{\delta |d_1|}  ] \nonumber\\
 &\leq& \ 2e^{-2} {\mathbb E} \Big[ \exp \Big( a G_{X_1}(X_1)\Big)\Big] \nonumber\\
 &\leq& \ 2e^{-2} K_1.
\end{eqnarray}
Combining the inequalities (\ref{finsa3f}) and (\ref{fisa3f}) together,  we obtain, for any $t \in [0, \delta )$,
\begin{eqnarray}
\mathbb{E}\,[e^{t d_1 } ] \ \leq\ 1+ \sum_{n=2}^{\infty} \frac{2}{e^{2}}\Big(\frac{t}{ \delta } \Big)^n  K_1
\ =\ 1+  \frac{2}{e^{2}} \frac{t^2 K_1 \delta^{-2} }{1-t \delta^{-1}}  \ \leq\  \exp \left( \frac{2}{e^{2}} \frac{t^2 K_1 \delta^{-2} }{1-t \delta^{-1}}   \right).
\end{eqnarray}
 Similarly, since $K_{n-i}(\rho)/K_{n-1}(\rho)\leq 1$ for all $i \in [2, n]$, we have,  for any $t \in [0,   \delta )$,
\begin{eqnarray}\label{fin36}
  \mathbb{E}\,[e^{t d_i  }|\mathcal{F}_{i-1}] \  \leq \ \exp \left(\frac{2}{e^{2}} \frac{t^2 K_2 \delta^{-2} }{1-t \delta^{-1}} \Big(\frac {K_{n-i}(\rho)}{K_{n-1}(\rho)}\Big)^2  \right) \, .
\end{eqnarray}
 By the tower property of conditional expectation, it follows that, for any $t \in [0, \delta )$,
\begin{eqnarray}
\mathbb{E}\,\big[e^{ tS_n} \big] &=&  \mathbb{E}\,\big[ \,\mathbb{E}\, [e^{ t S_n} |\mathcal{F}_{n-1}  ] \big] \nonumber\\
&=&  \mathbb{E}\,\big[ e^ { t S_{n-1}} \mathbb{E}\,  [e^ { t d_n}  |\mathcal{F}_{n-1}  ] \big] \nonumber \\
&\leq & \mathbb{E}\,\big[  e^ { t S_{n-1}}  \big] \exp \left( \frac{2}{e^{2}} \frac{t^2 K_2 \delta^{-2} }{1-t \delta^{-1}}  \right)  \nonumber\\
&\leq &  \exp \left( \frac{t^2 K  \delta^{-2} }{1-t \delta^{-1}}   \right),
\end{eqnarray}
where $$K=\frac{2}{e^{2}} \left( K_1+ K_2\sum_{i=2}^{n} \Big(\frac {K_{n-i}(\rho)}{K_{n-1}(\rho)}\Big)^2 \right).$$
Then using the exponential Markov  inequality, we deduce that,
for any $x\geq 0$ and $t \in [0, \delta )$,
\begin{eqnarray}
  \mathbb{P}\left( S_{n} \geq x \right)
 &\leq&  \mathbb{E}\, [e^{t\,(S_n -x) }  ] \nonumber\\
 &\leq&   \exp \left(-t x  +  \frac{t^2 K  \delta^{-2} }{1-t \delta^{-1}}   \right)\, .  \label{fisaxs}
\end{eqnarray}
The  minimum is reached at $$t=t(x):= \frac{x\delta^2/K}{x\delta/K + 1 + \sqrt{1+x\delta/K}} .$$ Substituting $t=t(x)$ in (\ref{fisaxs}),  we obtain the desired inequalities (\ref{Berie3}) and (\ref{Berie4}). \hfill\qed

\subsection{Qualitative results when ${\mathbb E} \big[ e^{a ( G_\varepsilon(\varepsilon))^p}\big]< \infty$ for $p>1$.}\label{QLW}
The next proposition follows easily from Theorem  3.2 of Liu and Watbled \cite{LW09}.
\begin{prop}\label{vprod}
Let $p>1$.
Assume that there exist some constants $a>0,$  $K_1\geq 1$ and $K_2\geq 1$ such that
\begin{equation}\label{laplacecd}
{\mathbb E} \Big[ \exp \Big( a \big(G_{X_1}(X_1)\big)^p\Big)\Big] \leq K_1
\quad \text{and} \quad
 {\mathbb E} \Big[ \exp \Big( a \big(G_\varepsilon(\varepsilon)\big)^p\Big)\Big] \leq K_2
 \, .
\end{equation}
 Let $q$ be the conjugate exponent of $p$ and let $\tau >0$ be such that
$$
(q \tau)^{\frac1 q} \left ( p a  \right)^{\frac 1 p}(1-\rho)=1 \, .
$$
Then, for any $\tau_1 > \tau$, there exist some positive numbers $t_1, x_1, A, B$,  depending only on $a, \rho,  K_1, K_2, p$ and $\tau_1$, such that
\begin{equation}\label{Bsesa2}
  \mathbb{E} \,[e^{\pm  tS_n}]\leq
  \begin{cases}
  \exp(n \tau_1 t^q) \quad  \text{ if\ \ $t \geq t_1$}\\
   \exp(n A t^2) \quad   \text{ if\ \ $t \in [0, t_1]$ }
  \end{cases}
\end{equation}
and
\begin{equation}\label{Ber2}
  {\mathbb P}\big(\pm  S_n \geq   x\big)\leq
  \begin{cases}
  \exp\left(-  a_1 x^p /n^{p-1}  \right) \quad  \text{if\ \ $x \geq nx_1$}\\
   \exp\left(- B x^2 /n \right) \quad   \text{if\ \ $x \in [0, nx_1]$,}
  \end{cases}
\end{equation}
where $a_1$ is such that $(q \tau)^{1/q} \left ( p a_1 \right)^{1/p}=1$.
\end{prop}
\begin{rem}
Assume that \eref{laplacecd} is satisfied for some $p\geq 1$.
From Proposition \ref{cram} (case $p=1$) and Proposition \ref{vprod} (case $p>1$), we infer that
 for any $x >0$, one can find a positive constant $c_x$  not depending on $n$ such  that
\begin{eqnarray}\label{jnks}
  \mathbb{P}\left(\pm   S_n \geq n x   \right)
   \leq   \exp\big( -c_x  n   \big) \, .
\end{eqnarray}
Moreover, for $x$ large enough, one can take $c_x= a_1 x^p$.
\end{rem}
\noindent\emph{Proof.} By condition (\ref{laplacecd}) and Proposition \ref{McD}, it follows that
\begin{equation}
 {\mathbb E} \Big[ \exp \Big( a (1-\rho)^p |d_1| ^p\Big)\Big] \leq  {\mathbb E} \Big[ \exp \Big( a \big( G_{X_1}(X_1)\big)^p\Big)\Big] \leq K_2  \nonumber
\end{equation}
and, for all $i\in [2, n], $
\begin{equation}
 {\mathbb E} \Big[ \exp \Big( a (1-\rho)^p | d_i |^p\Big) \ | \ \mathcal{F}_{i-1} \Big] \leq  {\mathbb E} \Big[ \exp \Big( a \big(  G_\varepsilon(\varepsilon)\big)^p \Big)\Big] \leq K_1  . \nonumber
\end{equation}
 Let $q> 1$ and $\tau >0$ be such that
$$
\frac1p + \frac1q=1\ \ \ \ \ \ \textrm{and} \ \ \ \  \ (q \tau)^{\frac1 q} \left ( p a  \right)^{\frac 1 p}(1-\rho)=1 \, .
$$
Then, by Theorem  3.2 of Liu and Watbled \cite{LW09}, for any $\tau_1 > \tau$, there exist $t_1, x_1, A, B>0$,  depending only on $a, \rho,  K_1, K_2, p$ and $\tau_1$, such that
  the claim of Proposition \ref{vprod} holds. \hfill\qed

\medskip

In particular, if $p=2$,  we have the following sub-Gaussian bound.
\begin{prop}\label{subgaussian}
Assume that there exist some constants $a>0,$ $K_1\geq 1$ and $K_2\geq 1$ such that
\begin{equation}\label{laplace2}
{\mathbb E} \Big[\exp \Big( a \big(G_{X_1}(X_1)\big)^2\Big)\Big] \leq K_1
\quad \text{and} \quad
 {\mathbb E} \Big[ \exp \Big( a \big( G_\varepsilon(\varepsilon)\big)^2\Big)\Big] \leq K_2
 \, .
\end{equation}
Then, there exists a constant $c>0$ depending only on $a, \rho, K_1$ and $K_2$ such that
\begin{equation}\label{Bdvbg}
  \mathbb{E} \,[e^{\pm  tS_n}] \leq \exp \big ( n\, c\, t^2) \quad \text{for all $t>0$.}
\end{equation}
Consequently, for any  $x>0$,
\begin{equation}\label{Besrs2}
  {\mathbb P}\big(\pm  S_n \geq  x\big)
  \leq \exp \left( -  \frac{x^2}{4nc}\right)\, .
\end{equation}
\end{prop}

\begin{rem}
As quoted at the beginning of Section \ref{deviationiq}, if $F$ satisfies only
\eref{contract}, Proposition \ref{subgaussian} holds provided \eref{weaker}
is satisfied with $f(x)=\exp(a x^2)$. This condition is implied by condition
\eref{DGWC},
which is due to Djellout {\it et al} \cite{DGW}.
\end{rem}
\noindent\emph{Proof.} Inequality (\ref{Bdvbg}) follows directly from (\ref{Bsesa2}).
Using the exponential Markov inequality, we deduce that
for any $x, t\geq 0$,
\begin{eqnarray}
  \mathbb{P}\left( S_{n} \geq x \right)
 &\leq&  \mathbb{E}\, [e^{t\,(S_n -x) }  ] \nonumber\\
 &\leq&   \exp \left(-t\, x  +  n\, c\,  t^2  \right)\, .  \label{figbxs}
\end{eqnarray}
The minimum is reached at $t=t(x):= x/(2 n  c)$. Substituting $t=t(x)$ in (\ref{figbxs}),  we obtain the desired inequality (\ref{Besrs2}). \hfill\qed

\subsection{Semi-exponential bound}
In the case where  $G_{X_1}(X_1)$ and $G_{\varepsilon}(\varepsilon)$
have semi-exponential moments, the following proposition holds. This proposition can be compared to the
corresponding results in Borovkov \cite{Bor} for partial sums of independent random variables,  Merlev\`{e}de {\it et al.}\ \cite{RPR10} for
partial sums of weakly
dependent sequences, and Fan {\it et al.}\ \cite{Fx1} for martingales.
\begin{prop}\label{findsa}
Let $p\in (0, 1)$.
Assume that there exist some positive constants $K_1$ and $K_2$ such that
\begin{equation}\label{laplace2}
{\mathbb E} \Big[ \big(G_{X_1}(X_1)\big)^2 \exp \Big( \big(G_{X_1}(X_1)\big)^p\Big)\Big] \leq K_1 \  \text{and} \ \
 {\mathbb E} \Big[ \big( G_\varepsilon(\varepsilon)\big)^2\exp  \Big(  \big( G_\varepsilon(\varepsilon)\big)^p\Big)\Big] \leq K_2
 \, .
\end{equation}
Set $$K= K_1+ K_2\sum_{i=2}^{n} \Big(\frac {K_{n-i}(\rho)}{K_{n-1}(\rho)}\Big)^2.$$
Then, for any $0\leq x<  K^{1/(2-p)}$,
\begin{eqnarray}\label{senib1}
  \mathbb{P}\left( \pm  S_n \geq x   \right)
 &\leq&    \exp\Bigg( -  \frac{x^2}{2K(K_{n-1}(\rho))^2}   \Bigg) \nonumber \\
 && + \big(K_{n-1}(\rho)\big)^{2}  \Big( \frac{ x^2 } {K^{1+p}}\Big)^{1/(1-p)}\exp\Bigg( -\Big(\frac{K}{x \big( K_{n-1}(\rho)\big)^{1-p}} \Big)^{p/(1-p)} \Bigg)\,  \
\end{eqnarray}
and, for any $x\geq K^{1/(2-p)}$,
\begin{eqnarray}\label{senib2}
  \mathbb{P}\left(\pm   S_n \geq x   \right)
 &\leq&    \exp\Bigg( -\bigg( \frac{x}{K_{n-1}(\rho)} \bigg)^{p} \bigg(1- \frac { K} {2 }\Big( \frac{K_{n-1}(\rho)}{ x} \Big)^{2-p}  \bigg)  \Bigg) \nonumber \\
 && + \ K \Big( \frac{ K_{n-1}(\rho) } {x}\Big)^{2}\exp\Bigg( -\Big(\frac{x}{K_{n-1}(\rho)} \Big)^{p} \Bigg)\, .
\end{eqnarray}
\end{prop}

\begin{rem}
In particular, there exists a positive constant $c$ such
that, for any $x >0$,
\begin{eqnarray}\label{jnks}
  \mathbb{P}\left(\pm   S_n \geq n x   \right)
   \leq C_x   \exp\big( -c\, x^p  n^p   \big),
\end{eqnarray}
where the constants  $C_x$ and $c$ do not depend on $n$.
\end{rem}

\begin{rem}
By a simple comparison, we find that for moderate
$ x \in (0,  K^{1/(2-p)})$, the second item in the right hand side of (\ref{senib1}) is less
than the first one. Thus for moderate $ x \in (0,  K^{1/(2-p)})$, the
 bound (\ref{senib1}) is a sub-Gaussian bound and is of the order
\begin{eqnarray}
  \exp\Bigg( -  \frac{x^2}{2K(K_{n-1}(\rho))^2}   \Bigg).
\end{eqnarray}
For  all $x\geq K^{1/(2-p)},$ bound (\ref{senib2}) is a semi-exponential bound and is of the order
\begin{eqnarray}\label{constant}
 \exp\Bigg( - \frac12 \bigg( \frac{x}{K_{n-1}(\rho)} \bigg)^{p}  \bigg)  \Bigg).
\end{eqnarray}
Moreover, when $x / K^{1/(2-p)} \rightarrow \infty,$ the constant $\frac12$ in (\ref{constant}) can be improved to $1+\varepsilon$ for any given $\varepsilon>0$.
\end{rem}

\noindent\emph{Proof.} The proof is  based on a truncation argument. For given $y>0,$
set $\eta_{i}=d_{i}\mathbf{1}_{\{d_{i}\leq y\}}$. Then $(\eta_{i}, \mathcal{F}_{i})_{i=1,...,n}$ is a sequence of supermartingale differences. Using a two term Taylor's expansion, we have, for all $t >0$,
\begin{equation*}
 e^{t \eta_i}  \leq 1+ t \eta_i +  \frac{t^2\eta_i^2}{2}\,  e^{ t  \eta_i}  \, .
\end{equation*}
Since $p \in (0,1)$, it follows that
$$
\eta_i^+= d_i{\bf 1}_{\{0 \leq d_i \leq y\}}\leq y \frac{d_i^p}{y^p}
{\bf 1}_{\{0 \leq d_i \leq y\}} \leq y^{1-p} (\eta_i^+)^p \, .
$$
Hence,
\begin{equation*}
 e^{t \eta_i}  \leq 1+ t \eta_i +  \frac{t^2\eta_i^2}{2}
  \exp\Big(  t y^{1-p} ( \eta_i ^+)^p \Big).
\end{equation*}
Since $\mathbb{E}[ \eta_i |\mathcal{F}_{i-1}]\leq \mathbb{E}[d_i |\mathcal{F}_{i-1}] = 0$, it follows that, for all $t >0$,
\begin{eqnarray*}
\mathbb{E}[e^{t \eta_i}|\mathcal{F}_{i-1}]  \leq  1+  \frac{t^2}{2} \mathbb{E}\Big[ \eta_i^2 \exp\Big(t  y^{1-p}  ( \eta_i ^+)^{p} \Big) \Big|\mathcal{F}_{i-1} \Big].
\end{eqnarray*}
By Proposition \ref{McD}, it follows that, for all $t>0$,
\begin{eqnarray*}
\mathbb{E}[e^{t \eta_1} ]  \leq  1+  \frac{t^2}{2} \mathbb{E}\Big[ \Big(K_{n-1}(\rho) G_{X_1}(X_1)   \Big)^2 \exp\Big(t  y^{1-p} \Big( K_{n-1}(\rho) G_{X_1}(X_1)   \Big)^{p} \Big)  \Big]
\end{eqnarray*}
and similarly, for $i \in [2, n],$
\begin{eqnarray*}
\mathbb{E}[e^{t \eta_i} |\mathcal{F}_{i-1}]  \leq  1+  \frac{t^2}{2} \mathbb{E}\Big[ \Big( K_{n-i}(\rho)  G_\varepsilon(\varepsilon)  \Big)^2 \exp\Big(t  y^{1-p} \Big( K_{n-i}(\rho)  G_\varepsilon(\varepsilon)  \Big)^{p} \Big)  \Big].
\end{eqnarray*}
Taking $ t=  y^{p-1}/\big( K_{n-1}(\rho) \big)^p$, by condition (\ref{laplace2}) and $K_{n-i}(\rho)/K_{n-1}(\rho) \leq 1$, we find that
\begin{eqnarray*}
\mathbb{E}[e^{t \eta_1} ]  &\leq&  1+  \frac{1 }{2}\Big(  \frac{y }{K_{n-1}(\rho) }    \Big)^{2p-2}  \mathbb{E}\Big[ \Big(   G_{X_1}(X_1)  \Big)^2 \exp\Big(  \big(  G_{X_1}(X_1)   \big)^{p} \Big)  \Big] \\
&\leq&  1+  \frac{ 1}{2}\Big(  \frac{y }{K_{n-1}(\rho) }    \Big)^{2p-2}   K_1 \\
&\leq& \exp\left( \frac{1}{2}  \Big(  \frac{y }{K_{n-1}(\rho) }    \Big)^{2p-2}  K_1 \right)
\end{eqnarray*}
and, for $i \in [2, n],$
\begin{eqnarray*}
\mathbb{E}[e^{t \eta_i} |\mathcal{F}_{i-1}]  &\leq&  1+ \frac{1 }{2}\Big(  \frac{y }{K_{n-1}(\rho) }    \Big)^{2p-2}  \Big(\frac {K_{n-i}(\rho)}{K_{n-1}(\rho)}\Big)^2  \mathbb{E}\Big[ \Big(     G_\varepsilon(\varepsilon)  \Big)^2 \exp\Big( \big(   G_\varepsilon(\varepsilon) \big)^{p} \Big)  \Big]\\
&\leq&  1+ \frac{1 }{2}\Big(  \frac{y }{K_{n-1}(\rho) }    \Big)^{2p-2}  K_2\Big(\frac {K_{n-i}(\rho)}{K_{n-1}(\rho)}\Big)^2 \ \\
&\leq& \exp\left( \frac{1 }{2}\Big(  \frac{y }{K_{n-1}(\rho) }    \Big)^{2p-2}   K_2\Big(\frac {K_{n-i}(\rho)}{K_{n-1}(\rho)}\Big)^2 \right).
\end{eqnarray*}
Hence, by the tower property of conditional expectation, it follows that
\begin{eqnarray}
\mathbb{E}\,\Big[e^{ t\sum_{i=1}^{n}\eta_i} \Big] &=&  \mathbb{E}\,\left[ \,\mathbb{E}\, [e^{ t\sum_{i=1}^{n}\eta_i} |\mathcal{F}_{n-1}  ] \right] \nonumber\\
&=&  \mathbb{E}\,\left[ e^ { t\sum_{i=1}^{n-1}\eta_i} \mathbb{E}\,  [e^ { t \eta_n}  |\mathcal{F}_{n-1}  ] \right] \nonumber \\
&\leq & \mathbb{E}\,\left[  e^ { t\sum_{i=1}^{n-1}\eta_i}  \right] \exp\left( \frac{1 }{2}\Big(  \frac{y }{K_{n-1}(\rho) }    \Big)^{2p-2}   K_2\Big(\frac {1}{K_{n-1}(\rho)}\Big)^2 \right) \nonumber\\
&\leq & \exp\left( \frac{1 }{2}\Big(  \frac{y }{K_{n-1}(\rho) }    \Big)^{2p-2} K   \right), \label{finew9}
\end{eqnarray}
where $$K= K_1+ K_2\sum_{i=2}^{n} \Big(\frac {K_{n-i}(\rho)}{K_{n-1}(\rho)}\Big)^2.$$
It is easy to see that
\begin{align}
  \mathbb{P}\left(  S_n \geq x   \right)
 &\leq  \mathbb{P}\left(   \sum_{i=1}^{n}\eta_i \geq x  \right)   +\  \mathbb{P}\left(   \sum_{i=1}^{n}d_{i}\mathbf{1}_{\{d_{i}> y\}} > 0  \right) \nonumber \\
 & \leq  \mathbb{P}\left(   \sum_{i=1}^{n}\eta_i \geq x  \right)   +
 \mathbb{P}\left( \max_{ 1\leq i\leq n}d_{i} > y \right)
=:  P_5 + \mathbb{P}\left( \max_{ 1\leq i\leq n}d_{i} > y \right). \label{fmuetoi}
\end{align}
For the first item of (\ref{fmuetoi}), by the  exponential Markov's inequality and (\ref{finew9}), we have
\begin{eqnarray}
P_5 &\leq&  \mathbb{E}\,\left[e^{ t(\sum_{i=1}^{n}\eta_i-x)}  \right]\ \ \leq \ \ \exp\left( -tx+  \frac{1 }{2}\Big(  \frac{y }{K_{n-1}(\rho) }    \Big)^{2p-2} K   \right). \label{fmusdv}
\end{eqnarray}
For the second item of (\ref{fmuetoi}),  we have the following estimation:
\begin{eqnarray*}
\mathbb{P}\left( \max_{ 1\leq i\leq n}d_{i} > y \right) &\leq& \sum_{i=1}^n \mathbb{P}\left(  d_{i} > y \right) \nonumber\\
&\leq& \sum_{i=1}^n \mathbb{P}\Big( \frac{d_{i}}{K_{n-1}(\rho) } > \frac{y}{K_{n-1}(\rho)}  \Big) \nonumber\\
 &\leq&  \frac{\exp\Big( -\big(y/K_{n-1}(\rho) \big)^{p} \Big) } {\big(y/K_{n-1}(\rho)\big)^{2} }  \ \sum_{i=1}^n \mathbb{E} \left[\Big(  \frac{d_{i}}{K_{n-1}(\rho)}\Big)^2 e^{ |d_{i} /K_{n-1}(\rho) |^{p}  } \right].
\end{eqnarray*}
By Proposition \ref{McD} and $K_{n-i}(\rho)/K_{n-1}(\rho) \leq 1$ again, it is easy to see that
\begin{multline*}
 \sum_{i=1}^n \mathbb{E} \left[\Big(  \frac{d_{i}}{K_{n-1}(\rho)}\Big)^2 e^{ |d_{i} /K_{n-1}(\rho) |^{p}  } \right] \\
\leq {\mathbb E} \Big[ \big(G_{X_1}(X_1)\big)^2 \exp \Big( \big(G_{X_1}(X_1)\big)^p\Big)\Big]  + \Big(\frac {K_{n-i}(\rho)}{K_{n-1}(\rho)}\Big)^2\sum_{i=2}^n {\mathbb E} \Big[ \big( G_\varepsilon(\varepsilon)\big)^2\exp  \Big(  \big( G_\varepsilon(\varepsilon)\big)^p\Big)\Big]
\leq K .
\end{multline*}
Thus
\begin{eqnarray}\label{fts90}
\mathbb{P}\left( \max_{ 1\leq i\leq n}d_{i} > y \right)
 &\leq& \frac{ K} {\big(y/K_{n-1}(\rho)\big)^{2} } \exp\Bigg( -\Big( \frac{y}{K_{n-1}(\rho)}  \Big)^{p} \Bigg) \, .
\end{eqnarray}
Combining (\ref{fmuetoi}), (\ref{fmusdv}) and (\ref{fts90}) together, it is to see that
\begin{eqnarray}
  \mathbb{P}\left(  S_n \geq x   \right)
 &\leq&  \exp\Bigg( -tx+  \frac{1 }{2}\Big(  \frac{y }{K_{n-1}(\rho) }    \Big)^{2p-2} K   \Bigg)   \nonumber\\
 && +\ \frac{ K} {\big(y/K_{n-1}(\rho)\big)^{2} } \exp\Bigg( -\Big(\frac{y}{K_{n-1}(\rho)}  \Big)^{p} \Bigg)\, .\ \   \label{finineq}
\end{eqnarray}
Recall that $t=  y^{p-1}/\big( K_{n-1}(\rho) \big)^p.$
Taking
\begin{eqnarray}  \label{f1}
y = \left\{ \begin{array}{ll}
\big(K/x\big)^{1/(1-p)} & \textrm{\ \ \ \ \ if $0\leq x < K^{1/(2-p)}$}, \\
x  & \textrm{\ \ \ \ \ if $x \geq K^{1/(2-p)}$, }
\end{array} \right.
\end{eqnarray}
we obtain the desired inequalities.  \hfill\qed


\subsection{McDiarmid inequality}\label{McDsec}
In this section, we consider the case where the increments $d_k$ are bounded.
We shall use an improved version of the well
known inequality by McDiarmid, which has been recently stated by Rio \cite{R13}.
For this inequality, we do not assume that \eref{c2} holds.
Hence, Proposition \ref{classic}  applies to any Markov chain $X_{n}=F(X_{n-1}, \varepsilon_n)$, for
$F$ satisfying \eref{contract}.

As in Rio \cite{R13},  let
$$\ell(t)= (t - \ln t -1) +t (e^t-1)^{-1} + \ln(1-e^{-t}) \ \ \ \ \
\textrm{for all}\ \ t > 0, $$ and let
$$\ell^*(x)=\sup_{t>0}\big(xt- \ell (t) \big) \ \ \ \ \ \textrm{for all}\ \ x > 0, $$ be the Young transform of $\ell(t)$.
As quoted by Rio \cite{R13}, the following inequalities hold
\begin{equation}\label{is50}
\ell^*(x) \geq (x^2-2x) \ln(1-x) \geq 2x^2 + x^4/6\, .
\end{equation}
Let also $(X_1', (\varepsilon'_i)_{i \geq 2})$ be an independent copy of
$(X_1, (\varepsilon_i)_{i \geq 2})$.
\begin{prop} \label{classic}
Assume that there exist some positive constants $M_k$ such that
\begin{equation}\label{in850}
\big \|d(X_1,X_1')
\big \|_\infty \leq M_1 \ \ \text{and} \ \
\big \|d\big(F(X_{k-1}, \varepsilon_k), F(X_{k-1}, \varepsilon'_k) \big)
\big \|_\infty \leq M_k \ \   \text{for $k \in [2,   n] $.}
\end{equation}
Let
$$
M^2(n,\rho)= \sum_{k=1}^n \big(K_{n-k}(\rho) M_k \big)^2 \quad \text{and} \quad D (n, \rho)=\sum_{k=1}^n K_{n-k}(\rho) M_k \, .
$$
Then, for any $t \geq 0$,
\begin{eqnarray}\label{rio1}
\mathbb{E}[e^{\pm   tS_n } ] \ \leq \  \exp \left (  \frac{D^2(n, \rho)}{M^2(n,\rho)}\  \ell \Big( \frac {M^2(n,\rho)\, x} {D(n, \rho)} \Big)\right) \,
\end{eqnarray}
and,  for any   $x \in [0, D(n, \rho)]$,
\begin{eqnarray}\label{rio2}
{\mathbb P}\big(\pm  S_n>x\big) \leq \exp \left ( -\frac{D^2(n, \rho)}{M^2(n,\rho)} \ \ell^*\Big( \frac {x} {D(n, \rho)} \Big)\right) \, .
\end{eqnarray}
Consequently, for any $x \in [0, D(n, \rho)]$,
\begin{eqnarray}\label{riosbound}
{\mathbb P}\big(\pm  S_n>x\big) &\leq&  \left ( \frac{D(n, \rho)-x}{D(n, \rho)}\right)^{ \frac{2D(n, \rho)x -x^2}{M^2(n,\rho)}} .
\end{eqnarray}
\end{prop}

\begin{rem} Since $(x^2-2x) \ln(1-x) \geq 2\,x^2$, inequality
(\ref{riosbound}) implies the following McDiarmid inequality
\begin{eqnarray*}
{\mathbb P}\big(\pm  S_n>x\big) \ \leq \  \exp \left ( -\frac{2 x^2}{M^2(n,\rho)}
 \right) \, .
\end{eqnarray*}
\end{rem}

\begin{rem}
Taking $\Delta (n, \rho)=K_{n-1}(\rho)\max_{1\leq k \leq n} M_k$, we obtain the upper bound: for any $x \in [0, n \Delta(n, \rho)]$,
\begin{eqnarray*}\label{rio2bis}
{\mathbb P}\big(\pm  S_n>x\big)
 \leq \exp \left ( -n \ell^*\Big( \frac {x} {n \Delta(n, \rho)} \Big)\right) \leq \exp \left ( -\frac{2 x^2}{n\Delta^2(n,\rho)}
 \right)
  \, .
\end{eqnarray*}
\end{rem}

\begin{rem}
If $F$ satisfies  \eref{c2}, then one can take
$M_1= \|d(X_1, X'_1)\|_\infty$ and $M_k=C \|\delta(\varepsilon_1, \varepsilon'_1)\|_\infty$ for $k \in [2, n]$.
\end{rem}

\medskip

\noindent\emph{Proof.}
Let
$$
u_{k-1}(x_1, \ldots, x_{k-1})= \text{ess}\inf_{\varepsilon_k} g_k(x_1, \ldots, F(x_{k-1}, \varepsilon_k))
$$
and
$$
v_{k-1}(x_1, \ldots, x_{k-1})= \text{ess}\sup_{\varepsilon_k} g_k(x_1, \ldots, F(x_{k-1}, \varepsilon_k))
$$
From the proof of Proposition \ref{McD},  it follows that
$$
u_{k-1}(X_1, \ldots, X_{k-1}) \leq d_k \leq v_{k-1}(X_1, \ldots, X_{k-1})\, .
$$
By Proposition \ref{McD} and  condition (\ref{in850}), we have
\[
v_{k-1}(X_1, \ldots, X_{k-1})-u_{k-1}(X_1, \ldots, X_{k-1}) \leq K_{n-k}(\rho) M_k\, .
\]
Now, following exactly the proof of  Theorem 3.1 of Rio \cite{R13} with $\Delta_k =  K_{n-k}(\rho) M_k$
we  obtain the inequalities (\ref{rio1}) and (\ref{rio2}).  Since
$\ell^*(x) \geq (x^2-2x) \ln(1-x)$, inequality (\ref{riosbound})
follows from (\ref{rio2}). \hfill\qed

\subsection{Fuk-Nagaev type bound}\label{Fuksec}
The next proposition follows easily from Corollary 2.3 of Fan {\it et al.} \cite{FGL12}.
\begin{prop}\label{po71}
Assume that  there exist two positive constants $V_1$ and  $V_2$ such that
$$
{\mathbb E} \big[  \big(  G_{X_1}(X_1)\big)^2\big] \leq  V_1 \quad \text{and} \quad
{\mathbb E} \big[  \big(   G_\varepsilon(\varepsilon)\big)^2\big] \leq
  V_2  \, .
 $$
Let
\begin{equation}\label{hoeffineq}
V=  V_1\Big(K_{n-1}(\rho)\Big)^2 +  V_2 \sum_{i=2}^n \Big(K_{n-i}(\rho)\Big)^2.
\end{equation}
Then, for any $x,y > 0$,
\begin{multline}\label{Hineq}
{\mathbb P}\big(\pm  S_n>x \big)  \leq    H_n \left(\frac{x }{y K_{n-1}(\rho)} , \frac{\sqrt{V} }{y K_{n-1}(\rho) }\right)
   + {\mathbb P}\left( \max\left \{ G_{X_1}(X_1), \max_{2 \leq i \leq n} G_\varepsilon(\varepsilon_i)\right \} > y \right) \,  ,
\end{multline}
where
\begin{eqnarray}\label{fkmns}
H_{n}(x,v)=\left\{\left( \frac{v^2}{
x+v^2}\right)^{ x+v^2 }\left( \frac{n}{n-x}\right)^{ n-x }
\right\}^{\frac{n}{n+v^2} }\mathbf{1}_{\{x \leq n\}}
\end{eqnarray}
with the convention that $(+\infty)^0=1$   \emph{(}which  applies when $x=n$\emph{)}.
\end{prop}

\noindent\emph{Proof.} We apply Corollary 2.3 of Fan {\it et al.} \cite{FGL12} with the truncature level $y K_{n-1}(\rho)$.
By Proposition \ref{McD}, $|d_1| \leq K_{n-1}(\rho)G_{X_1}(X_1)$ and $|d_i| \leq K_{n-i}(\rho) G_\varepsilon(\varepsilon_i)$ for $i\in  [2, n]$. Hence
$$
{\mathbb E}\big[d_1^2 {\bf 1}_{\{d_1 \leq y K_{n-1}(\rho)\}}\big]
  \leq  \Big(K_{n-1}(\rho)\Big)^2 {\mathbb E} \big[  \big(  G_{X_1}(X_1)\big)^2\big] \leq  \Big(K_{n-1}(\rho)\Big)^2 V_1 \,
$$
and, for $i\in [2, n]$,
$$
 {\mathbb E}\big[d_i^2 {\bf 1}_{\{d_i \leq y K_{n-1}(\rho)\}}|{\mathcal F}_{i-1}\big]
  \leq \Big(K_{n-i}(\rho)\Big)^2
 {\mathbb E} \big[  \big( G_\varepsilon(\varepsilon)\big)^2\big] \leq
  \Big(K_{n-i}(\rho)\Big)^2 V_2\, .
$$
It follows from  Corollary 2.3 of Fan  {\it et al.} \cite{FGL12} that
$$
{\mathbb P}(S_n>x)  \leq    H_n \left(\frac{x }{y K_{n-1}(\rho)} , \frac{\sqrt{V} }{y K_{n-1}(\rho) }\right) \\
+ {\mathbb P}\bigg( \max_{1 \leq i \leq n} d_i >  y K_{n-1}(\rho) \bigg)\, .
$$
Inequality \eref{Hineq} follows by applying
Proposition \ref{McD} again. \hfill \qed

\medskip

In particular, if $G_{X_1}(X_1)$ and $G_\varepsilon(\varepsilon)$ are bounded, then Proposition \ref{po71} implies the following Hoeffding bound.
\begin{prop}\label{fhoedding}
Assume that there exist some positive constants $M$, $V_1$ and $V_2$ such that
\begin{equation*}\label{conditionHoeffding}
 G_{X_1}(X_1) \leq M, \quad G_\varepsilon(\varepsilon) \leq M , \quad
 {\mathbb E}\big[\big(G_{X_1}(X_1)\big)^2 \big] \leq V_1 \quad \text{and} \quad
 {\mathbb E}\big[\big(G_\varepsilon(\varepsilon)\big)^2 \big] \leq V_2
 \, .
\end{equation*}
  Then, for any $x > 0$,
\begin{eqnarray}\label{H}
{\mathbb P}\big(\pm S_n>x\big)  &\leq&    H_n \left(\frac{x }{M K_{n-1}(\rho)} , \frac{\sqrt{V} }{M K_{n-1}(\rho) }\right) ,
\end{eqnarray}
where $H_{n}(x,v)$ and $V$ are defined by (\ref{fkmns}) and (\ref{hoeffineq}), respectively.
\end{prop}

\begin{rem}
According to Remark 2.1 of Fan  {\it et al.} \cite{FGL12},  for any $x\geq0$ and any $v>0$, it holds
\begin{eqnarray}
H_{n}(x, v)&\leq& B(x,v):=\left(\frac{v^2}{x+v^2} \right)^{x+v^2}e^x  \label{freedma1} \\
&\leq& B_1(x,v):= \exp\left\{-\frac{x^2}{2(v^2+\frac{1}{3}x )}\right\} . \label{Bernstein}
\end{eqnarray}
Note that (\ref{freedma1}) and (\ref{Bernstein}) are respectively known as   Bennett's
and Bernstein's bounds.
Then, inequality (\ref{H}) also implies  Bennett's and Bernstein's bounds
\begin{eqnarray*}
\mathbb{P}\big(\pm S_n>x\big) \  \leq  \  B\left(\frac{x }{M K_{n-1}(\rho)} , \frac{\sqrt{V} }{M K_{n-1}(\rho) }\right)\  \leq \   B_1\left(\frac{x }{M K_{n-1}(\rho)} , \frac{\sqrt{V} }{M K_{n-1}(\rho) }\right).
\end{eqnarray*}
\end{rem}

\medskip

We now consider the case where the random variables
$G_{X_1}(X_1)$ and $G_\varepsilon(\varepsilon)$ have only a weak moment of order $p > 2$.
For any real-valued random variable $Z$ and any $p\geq 1$, define the weak
moment of order $p$  by
\begin{equation}\label{weakp}
\|Z\|_{w,p}^p=\sup_{x>0} x^p{\mathbb P}(|Z|>x)\, .
\end{equation}
\begin{prop}\label{fnhbslp}
 Let $p > 2$.
Assume that there exist some positive constants  $V_1,$ $V_2,$
$A_1(p)$ and $A_2(p)$ such that
\begin{eqnarray*}
&&{\mathbb E}\big[\big(G_{X_1}(X_1)\big)^2 \big] \leq V_1 \ , \ \ \ \ \
\quad  \quad \ \ \ \ \
 {\mathbb E}\big[\big( G_\varepsilon(\varepsilon)\big)^2 \big] \leq V_2
\, , \\
&&\Big \|   G_{X_1}(X_1)\Big \|_{w,p}^p \leq A_1(p)\,   \
 \quad  \ \text{and} \ \   \quad
  \Big \| G_\varepsilon(\varepsilon)\Big\|_{w,p}^p \leq A_2(p) \, .
\end{eqnarray*}
Let $V$ be defined by (\ref{hoeffineq}),
and let
 \begin{equation*}
A(p)= A_1(p) +
(n-1)A_2(p) \, .
\end{equation*}   Then, for any $x,y > 0$,
\begin{eqnarray}\label{weakpFN}
{\mathbb P}\big(\pm  S_n>x\big)  &\leq&    H_n \left(\frac{x}{yK_{n-1}(\rho)} , \frac{\sqrt{V} }{y K_{n-1}(\rho)}\right)    +  \frac{A(p)}{y^p },
\end{eqnarray}
where $H_{n}(x,v)$ is defined by (\ref{fkmns}).
\end{prop}
\begin{rem}
Assume that $G_{X_1}(X_1)$ and $G_\varepsilon(\varepsilon)$ have a weak moment of order
$p>2$.
Taking
$$
y= \frac{3 n x}{2pK_{n-1}(\rho)\ln(n)}
$$
in inequality \eref{weakpFN},
we infer that, for any $x>0$,
$$
{\mathbb P}\big( \pm S_n>nx \big) \leq \frac{C_x (\ln(n))^p}{n^{p-1}} \, ,
$$
for some positive $C_x$ not depending on $n$.
\end{rem}

\medskip

If the martingale differences $d_i$ have $p$th moments ($p\geq 2$), then we have the following Fuk-type inequality (cf.\ Corollary $3'$ of Fuk \cite{F73}).
\begin{prop}
 Let $p\geq 2$. Assume that there exist some positive constants   $V_1,$ $V_2,$
$A_1(p)$ and $A_2(p)$ such that
\begin{eqnarray}
&&{\mathbb E}\big[\big(G_{X_1}(X_1)\big)^2 \big] \leq V_1\, ,
\ \ \ \ \ \ \ \ \quad   \ \ \ \ \
{\mathbb E}\big[\big(G_\varepsilon(\varepsilon)\big)^2 \big] \leq V_2\, ,
 \, \nonumber \\
&&{\mathbb E}\big[\big(G_{X_1}(X_1)\big)^p \big] \leq A_1(p)\,
\quad   \text{and}  \quad \ \ \,
 {\mathbb E}\big[\big(G_\varepsilon(\varepsilon)\big)^p \big] \leq A_2(p)
  \,  .\label{sdvcx}
\end{eqnarray}
Let $V$ be defined by \eref{hoeffineq}, and let  $$A(p)=  A_1(p)\big(K_{n-1}(\rho) \big)^p   + A_2(p)\sum_{i=2}^n\big(K_{n-i}(\rho) \big)^p. $$  Then, for any $x  > 0$,
\begin{eqnarray}\label{fuki}
{\mathbb P}( |S_n| >x)  &\leq& 2\Big(1+ \frac2 p \Big)^p \frac{A(p)}{x^p}   + 2 \exp\left( - \frac{2}{(p+2)^2 e^p}  \frac{x^2}{V}  \right).
\end{eqnarray}
\end{prop}
\begin{rem}
Since  $A(p)$ is of order $n$,  it easy to see that
the term
$$\exp\left( - \frac{2}{(p+2)^2 e^p}  \frac{(nx)^2}{V}  \right)$$ is decreasing at an exponential order, and
that  the term $$2\Big(1+ \frac2 p \Big)^p \frac{A(p)}{(xn)^p}$$
 is of order $n^{1-p}$.
 Thus, for any $x>0$ and all $n$,
$$
  {\mathbb P}(|S_n|>nx)\leq \frac{C_x}{ n^{p-1}} \, ,
$$
for some positive $C_x$ not depending on $n$. Note that the last inequality is optimal under the stated condition, even if $S_n$ is a sum of iid random variables.
\end{rem}

\noindent\emph{Proof.}
By Proposition \ref{McD}  and condition (\ref{sdvcx}), it follows that
\begin{eqnarray*}
\sum_{i=1}^n\mathbb{E}[|d_i|^p | \mathcal{F}_{i-1} ] &\leq& \mathbb{E}[| K_{n-1}(\rho) G_{X_1}(X_1)|^p ] + \sum_{i=2}^n\mathbb{E}[| K_{n-i}(\rho) G_{\varepsilon}(\varepsilon_i)|^p ] \\
&\leq & \big(K_{n-1}(\rho) \big)^p \mathbb{E}[|G_{X_1}(X_1)|^p ] + \sum_{i=2}^n\big(K_{n-i}(\rho) \big)^p\mathbb{E}[| G_{\varepsilon}(\varepsilon_i)|^p ]\\
&\leq& A_1(p)\big(K_{n-1}(\rho) \big)^p   + A_2(p)\sum_{i=2}^n\big(K_{n-i}(\rho) \big)^p
=A(p).
\end{eqnarray*}
Notice that $A(2)=V$. Using Corollary $3'$ of Fuk \cite{F73}, we obtain the desired inequality. \hfill\qed

\subsection{von Bahr-Esseen bound}\label{VBEB}
In the first proposition of this section, we assume that the dominating random variables $G_{X_1}(X_1)$
and
$G_\varepsilon(\varepsilon_k)$ have only a moment of order $p \in [1,2]$.
For similar inequalities in the case where the
$X_i$'s are independent, we refer to Pinelis \cite{P10}.
\begin{prop}\label{VBEI}
Let $p \in [1, 2]$. Assume that
\begin{equation}\label{fnf58}
{\mathbb E} \Big[  \Big(  G_{X_1}(X_1)\Big)^p\Big] \leq A_1(p)
\quad  \text{and} \quad
{\mathbb E} \Big[  \Big(   G_\varepsilon(\varepsilon)\Big)^p\Big] \leq A_2(p) \, .
\end{equation}
Then
\begin{equation}\label{vBE}
  \| S_n \|_p  \leq    \big( A(n,\rho,p) \big)^{1/p},
\end{equation}
 where
\begin{equation}\label{Ap}
A(n,\rho,p)= A_1(p)\big(K_{n-1}(\rho)\big)^p    +
2^{2-p} A_2(p)  \sum_{k=2}^n \big(K_{n-k}(\rho)\big)^p \, .
\end{equation}
\end{prop}
\begin{rem}
The constant $2^{2-p}$ in \eref{Ap} can be replaced by the more precise constant $\tilde C_p$ described in Proposition 1.8 of Pinelis \cite{P10}.
\end{rem}
\begin{rem}
Assume that  $F$ satisfies only \eref{contract}. Then, it follows from the proof of Proposition \ref{VBEI} that
the inequality \eref{vBE} remains true if the second condition of  \eref{fnf58} is replaced by
\begin{equation*}
\sup_{k \in [2,n]}
  {\mathbb E} \Big [ \Big ( H_\varepsilon(X_{k-1}, \varepsilon_k) \Big)^p
 \Big ] \leq A_2(p) \, .
\end{equation*}
\end{rem}
\noindent\emph{Proof.} Using an improvement of the von Bahr-Esseen inequality  (see  inequality (1.11) in Pinelis \cite{P10}), we have
$$
\| S_n \|_p^p \leq \|d_1\|_p^p +  \tilde C_p \sum_{k=2}^n  \|\, d_k \|_p^p \, ,
$$
where the constant $\tilde C_p$ is described in Proposition 1.8 of Pinelis \cite{P10}, and is such that $\tilde C_p \leq 2^{2-p}$
for any $p \in [1,2]$.
By Proposition \ref{McD}, it follows that
\begin{eqnarray*}
\| S_n \|_p^p   &\leq&  \, \bigg( \,  \big(K_{n-1}(\rho)\big)^p   {\mathbb E} \Big[  \Big(  G_{X_1}(X_1)\Big)^p\Big]    + \tilde C_p
    \sum_{k=2}^n \big(K_{n-k}(\rho)\big)^p {\mathbb E} \Big[  \Big(   G_\varepsilon(\varepsilon)\Big)^p\Big]\, \bigg) \\
    &\leq&  \bigg( A_1(p)\big(K_{n-1}(\rho)\big)^p    +
\tilde C_p  A_2(p)  \sum_{k=2}^n \big(K_{n-k}(\rho)\big)^p \, \bigg),
\end{eqnarray*}
which gives the desired inequality. \hfill\qed

\medskip

We now consider the case where the variables $G_{X_1}(X_1)$ and
$G_\varepsilon(\varepsilon_k)$ have only a weak moment of order $p \in (1,2)$.
Recall that the weak moment $\| Z \|_{w,p}^p$ has been defined by \eref{weakp}.

\begin{prop}\label{weakVBEI}
Let $p \in (1, 2)$. Assume that
\begin{equation}\label{weakVB}
\Big \|   G_{X_1}(X_1)\Big \|_{w,p}^p \leq A_1(p)
 \quad  \text{and} \quad
  \Big \|G_\varepsilon(\varepsilon)\Big\|_{w,p}^p \leq A_2(p)
  \, .
\end{equation}
Then, for any $x>0,$
\begin{equation}\label{weakVBE}
 {\mathbb P}(|S_n|>x) \leq  \frac{C_p B(n,\rho,p)}{x^p}\, ,
\end{equation}
where
 $$
 C_p= \frac{4p}{(p-1)} +\frac{8p}{(p-2)}\,
$$
and
 $$
B(n,\rho,p)= A_1(p)\big(K_{n-1}(\rho)\big)^p    +
 A_2(p)  \sum_{k=2}^n \big(K_{n-k}(\rho)\big)^p  \,.
$$
\end{prop}

\begin{rem}
Assume that  $F$ satisfies only \eref{contract}. Then, it follows from the proof of Proposition \ref{weakVBEI} that
the inequality
 \eref{weakVBE} remains true if the second condition of  \eref{weakVB} is replaced by
\begin{equation*}
\sup_{k \in [2,n]}
  \Big\| H_\varepsilon(X_{k-1}, \varepsilon_k) \Big\|_{w,p}^p \leq A_2(p) \, .
\end{equation*}
\end{rem}

\noindent\emph{Proof.}
This proof is based on a truncation argument. For given $x>0,$ let
\begin{eqnarray*}
&& \xi_1 = d_1 \textbf{1}_{\{d_1 \leq x\}}\, , \, \ \ \  \ \ \ \  \  \  \ \ \, \xi_1' = d_1 \textbf{1}_{\{d_1 > x \}}\, , \\
&& \xi_k = d_k \textbf{1}_{\{d_k\leq x \}} \ \  \ \   \textrm{and}  \ \  \ \ \ \xi_k' = d_k \textbf{1}_{\{d_k>  x\}}\,.
\end{eqnarray*}
Define
\begin{eqnarray*}
&&\eta_1= \xi_1-{\mathbb E}[\xi_1]\, ,  \ \ \  \ \ \ \  \  \  \ \ \ \ \ \ \ \ \ \ \  \eta'_1= \xi'_1-{\mathbb E}[\xi'_1]\, , \\
&&\eta_k= \xi_k-{\mathbb E}[\xi_k|{\mathcal F}_{k-1}] \ \  \ \   \textrm{and}  \ \  \ \ \ \eta'_k= \xi'_k-{\mathbb E}[\xi'_k|{\mathcal F}_{k-1}]\,.
\end{eqnarray*}
It is obvious that
\begin{equation}\label{firstbound}
{\mathbb P}(|S_n|>x) \leq {\mathbb P} \Big(  \Big |\sum_{k=1}^n \eta_k \Big | > \frac x
 2\Big)
 +
 {\mathbb P} \Big(  \Big |\sum_{k=1}^n \eta'_k \Big | > \frac x
 2\Big).
\end{equation}
Applying Markov's inequality, we get
\begin{equation}\label{M1}
{\mathbb P} \Big(  \Big |\sum_{k=1}^n \eta'_k \Big | > \frac x
 2\Big) \leq \frac{2}{x} \sum_{k=1}^n \|\eta'_k\|_1 \leq \frac{4}{x} \sum_{k=1}^n \|\xi'_k\|_1 \, .
\end{equation}
Recall that, if $Z$ is any real-valued random variable such that
\begin{equation}\label{quant}
{\mathbb P}(|Z|>x) \leq H(x)
\end{equation}
for a tail function $H$, then
\begin{equation}\label{basic1}
{\mathbb E}(|Z|\textbf{1}_{\{|Z|>a\}}) \leq \int_0^{H(a)} Q(u) du \, ,
\end{equation}
where $Q$ is the cadlag inverse of $H$.
Using Proposition \ref{McD}, we have
\begin{equation}\label{domin}
{\mathbb P}(|d_k|>x) \leq H_k(x),
\end{equation}
where $H_1(x)= \min \{1, x^{-p} A_1(p) (K_{n-1}(\rho))^p\}$ and
$H_k(x)= \min \{1, x^{-p} A_2(p) (K_{n-k}(\rho))^p\}$ if $k\in [2, n]$.
Hence, applying \eref{basic1},  we obtain
\begin{equation}\label{1.1}
\|\xi'_1\|_1 \leq (A_1(p))^{1/p} K_{n-1}(\rho)
\int_0^{H_1(x)} u^{-1/p}  du \leq \frac{p}{p-1} A_1(p) (K_{n-1}(\rho))^p x^{1-p} \, .
\end{equation}
Similarly, for $k \in [2,n]$,
\begin{equation}\label{1.k}
\|\xi'_k\|_1 \leq (A_2(p))^{1/p} K_{n-k}(\rho)
\int_0^{H_k(x)} u^{-1/p}  du \leq \frac{p}{p-1} A_2(p) (K_{n-k}(\rho))^p x^{1-p} \, .
\end{equation}
Consequently, from \eref{M1}, \eref{1.1} and \eref{1.k},
\begin{equation}\label{firststep}
{\mathbb P} \Big(  \Big |\sum_{k=1}^n \eta'_k \Big | > \frac x
 2\Big) \leq \frac{4p B(n, \rho, p)}{(p-1)x^p} \, .
\end{equation}
On the other hand, the $\eta_k$'s being martingales differences,
\begin{equation}\label{M2}
{\mathbb P} \Big(  \Big |\sum_{k=1}^n \eta_k \Big | > \frac x
 2\Big) \leq \frac{4}{x^2} \sum_{k=1}^n \|\eta_k\|_2^2
 \leq \frac{4}{x^2} \sum_{k=1}^n \|\xi_k\|_2^2
\end{equation}
Recall that, if $Z$ is any real-valued random variable satisfying \eref{quant},
\begin{equation}\label{basic2}
{\mathbb E}(Z^2 {\bf 1}_{|Z|\leq a})\leq
{\mathbb E}((Z\wedge a)^2) \leq \int_0^1 \min \{ Q^2(u), a^2 \} du
\leq 2 \int_{H(a)}^1 Q^2(u) du \, .
\end{equation}
Using \eref{domin} and \eref{basic2},
we obtain
\begin{equation}\label{2.1}
\|\xi_1\|_2^2 \leq 2(A_1(p))^{2/p} (K_{n-1}(\rho))^2 \int_{H_1(x)}^1 u^{-2/p} du
\leq \frac{2p}{2-p} A_1(p) (K_{n-1}(\rho))^p x^{2-p} \, .
\end{equation}
Similarly, for $k \in [2,  n]$,
\begin{equation}\label{2.k}
\|\xi_k\|_2^2 \leq 2(A_2(p))^{2/p} (K_{n-k}(\rho))^2 \int_{H_k(x)}^1 u^{-2/p} du
\leq \frac{2p}{2-p} A_2(p) (K_{n-k}(\rho))^p x^{2-p}\, .
\end{equation}
Consequently, from \eref{M2}, \eref{2.1} and \eref{2.k},
\begin{equation}\label{secondstep}
{\mathbb P} \Big(  \Big |\sum_{k=1}^n \eta_k \Big | > \frac x
 2\Big) \leq \frac{8p B(n, \rho, p)}{(p-2)x^p} \, .
\end{equation}
Inequality \eref{weakVBE} follows from
\eref{firstbound},  \eref{firststep} and \eref{secondstep}. \qed

\subsection{Marcinkiewicz-Zygmund bound}\label{MZB}
We now  assume that the dominating random variables $G_{X_1}(X_1)$
and
$G_\varepsilon(\varepsilon_k)$ have a moment of order $p \geq 2$.
\begin{prop}\label{MZIP}
Let $p \geq 2$. Assume that
\begin{equation}\label{shnbs}
{\mathbb E} \Big[  \Big(  G_{X_1}(X_1)\Big)^p\Big] \leq A_1(p)
 \quad  \text{and} \quad
{\mathbb E} \Big[  \Big(   G_\varepsilon(\varepsilon)\Big)^p\Big] \leq A_2(p)\, .
\end{equation}
Then
\begin{equation}\label{MZI}
  \| S_n \|_p \leq   \sqrt{ A(n,\rho,p)}\, ,
\end{equation}
 where
 $$
A(n,\rho,p)= \big(K_{n-1}(\rho) \big)^2 \big(A_1(p) \big)^{2/p} +
(p-1)\big(A_1(p) \big)^{2/p} \sum_{k=2}^n \big(K_{n-k}(\rho) \big)^2 \, .
$$
\end{prop}
\begin{rem}
Assume that  $F$ satisfies only \eref{contract}. Then, it follows from the proof of Proposition \ref{MZIP} that
the inequality
\eref{MZI} remains true if the second condition of  \eref{shnbs} is replaced by
\begin{equation*}
\sup_{k \in [2,n]}
  {\mathbb E} \Big [ \Big ( H_\varepsilon(X_{k-1}, \varepsilon_k) \Big)^p
 \Big ] \leq A_2(p) \, .
\end{equation*}
\end{rem}
\noindent\emph{Proof.} Using Theorem 2.1 of Rio \cite{R09}, we have
 $$
\| S_n \|_p^2  \leq \|d_1\|_p^2 +(p-1) \sum_{k=2}^n  \| d_k \|_p^2 \, .
$$
By Proposition \ref{McD}  and condition (\ref{shnbs}), it follows that
\begin{eqnarray*}
\| S_n \|_p^2   &\leq&   \big(K_{n-1}(\rho)\big)^2  \left( {\mathbb E} \Big[  \Big(  G_{X_1}(X_1)\Big)^p\Big]  \right)^{2/p} + (p-1)
    \sum_{k=2}^n \big(K_{n-k}(\rho)\big)^2 \left({\mathbb E} \Big[  \Big(   G_\varepsilon(\varepsilon)\Big)^p\Big]\right)^{2/p} \\
    &\leq&  A(n,\rho,p),
\end{eqnarray*}
which gives the desired inequality. \hfill\qed

\subsection{Burkholder-Rosenthal bounds}\label{BRB}
When  the dominating random variables $G_{X_1}(X_1)$
and
$G_\varepsilon(\varepsilon_k)$ have  a moment of order $p \geq 2$,
one can prove the following proposition. For similar inequalities in the case where the
$X_i$'s are independent, we refer to Pinelis \cite{P13}.
\begin{prop}\label{Ros}
Assume that there exist two constants $ V_1 \geq 0$ and $V_2 \geq 0$ such that
\begin{equation}\label{Strongros}
{\mathbb E}\big[\big(G_{X_1}(X_1)\big)^2 \big] \leq V_1
\quad \text{and} \quad
 {\mathbb E}\big[\big( G_\varepsilon(\varepsilon)\big)^2 \big] \leq V_2
\, .
\end{equation}
Let
\begin{equation}\label{variance}
V= V_1 \big(K_{n-1}(\rho) \big)^2  +
V_2 \sum_{k=2}^n \big(K_{n-k}(\rho) \big)^2 \, .
\end{equation}
For any  $p\geq 2$, there exist two positive constants $C_1(p)$ and $C_2(p)$ such that
\begin{equation}\label{rosI}
  \| S_n \|_p \leq   C_1(p) \sqrt{ V} + C_2(p)  \Big \| \max\Big \{ K_{n-1}(\rho)G_{X_1}(X_1), \max_{2 \leq i \leq n} K_{n-i}(\rho)G_\varepsilon(\varepsilon_i)\Big \}\Big \|_p  .
\end{equation}
\end{prop}

\begin{rem}
According to the proof of Theorem 4.1 of Pinelis \cite{P94}, one can take $C_1(p)=60c$ and $C_2(p)= 120\sqrt{c} e^{p/c}$ for any $c \in [1,p]$.
\end{rem}

\medskip

\noindent\emph{Proof.}
Applying Proposition \ref{McD}, we have
$|d_1| \leq K_{n-1}(\rho) G_{X_1}(X_1)$ and $|d_k| \leq K_{n-k}(\rho) G_{\varepsilon}(\varepsilon_k)$ for $k \in [2,  n],$
and consequently
$$
{\mathbb E}[d_1^2] \leq \big(K_{n-1}(\rho) \big)^2 V_1 \quad \text{and} \quad {\mathbb E}[d_k^2|{\mathcal F}_{k-1}] \leq  \big(K_{n-k}(\rho)\big)^2 V_2\, \text{ for $k\in [2,  n].$}
$$
Then the proposition  follows  directly from Theorem 4.1 of Pinelis \cite{P94}. \hfill \qed

\medskip

We now consider the case where the random variables
$G_{X_1}(X_1)$ and $G_\varepsilon(\varepsilon)$ have  a weak moment of order $p >2$.
Recall that the weak moment $\| Z \|_{w,p}^p$ has been defined by \eref{weakp}.
\begin{prop}\label{RCG}
Assume that \eref{Strongros} holds, and let $V$ be defined by \eref{variance}.
Then, for any  $p\geq 2$, there exist two positive constants $C_1(p)$ and $C_2(p)$ such that
\begin{equation}\label{weakrosI}
  {\mathbb P}(|S_n|>t) \leq
  \frac{1}{t^p} \left \{  C_1(p) {V}^{p/2} + C_2(p)  \Big \| \max\Big \{ K_{n-1}(\rho)G_{X_1}(X_1), \max_{2 \leq i \leq n} K_{n-i}(\rho)G_\varepsilon(\varepsilon_i)\Big \}\Big \|^p_{w,p} \right \} .
\end{equation}
\end{prop}

\begin{rem}
Assume that  $F$ satisfies only \eref{contract}. Then, it follows from the proofs of Propositions \ref{Ros} and \ref{RCG} that
the  inequalities \eref{rosI}  and  \eref{weakrosI} remain true if the second condition of  \eref{Strongros} is replaced by
\begin{equation*}
\sup_{k \in [2,n]}
   \Big \| {\mathbb E}\Big [\big( H_\varepsilon(X_{k-1}, \varepsilon_k)\big)^2\Big | X_{k-1} \Big ]\Big \|_{p/2}
  \leq V_2 \, ,
\end{equation*}
and by taking $H_\varepsilon(X_{k-1}, \varepsilon_k)$ instead of $G_\varepsilon(\varepsilon_k)$ in the second terms on right hand
of \eref{rosI}  and  \eref{weakrosI}.
\end{rem}

\noindent\emph{Proof.} It is the same as that of Proposition \ref{Ros}, by applying Theorem 6.3 in Chazottes and Gou\"ezel \cite{CG12}. \hfill \qed

\section{Application to the Wasserstein distance between the empirical distribution
and the invariant distribution}
\setcounter{equation}{0}
\subsection{Definition and  upper bounds}\label{DU}
Recall that the Wasserstein distance $W_1(\nu_1, \nu_2)$ between two probability measures $\nu_1, \nu_2$
on $({\mathcal X},d)$ is defined by
$$
  W_1(\nu_1, \nu_2) = \inf_{\lambda \in M(\nu_1, \nu_2)} \int d(x,y) \lambda(dx, dy) \, ,
$$
where $M(\nu_1, \nu_2)$ is the set of probability measures on
${\mathcal X}\times {\mathcal X}$ with margins $\nu_1$ and $\nu_2$.

Let $\Lambda_1({\mathcal X})$ be the set of functions from $({\mathcal X},d)$ to
${\mathbb R}$ such that
$$
|g(x)-g(y)| \leq d(x,y) \, .
$$
Recall that $W_1(\nu_1, \nu_2)$ can be expressed  {\it via} its dual form (see for instance
the equality (5.11) in Villani \cite{V})
$$
  W_1(\nu_1, \nu_2) = \sup_{g \in \Lambda_1({\mathcal X})} |\nu_1(g)- \nu_2(g)| \, .
$$

Let $\mu_n$ be the empirical distribution of the random variables $X_1, X_2, ..., X_n$, that
is
$$
\mu_n= \frac 1 n \sum_{k=1}^n \delta_{X_k} \, ,
$$
and let $\mu$ be the unique invariant distribution of the chain. It is easy to see
that the function $f$ defined by
$$
nW_1(\mu_n, \mu)= f(X_1, X_2, \ldots, X_n) := \sup_{g \in \Lambda_1({\mathcal X})} \,\Big|\sum_{i=1}^n \big(g(X_i)-\mu(g)\big)\Big|\, ,
$$
is separately Lipschitz, and satisfies \eref{codiMD}. Hence, all the inequalities of
Section \ref{deviationiq} apply to
$$
S_n= nW_1(\mu_n, \mu)- n{\mathbb E}[W_1(\mu_n, \mu)] \, .
$$
Let us only give some qualitative consequences of these inequalities:
\begin{itemize}
\item If \eref{laplacecd} holds for some $p\geq 1$, then
there exist some positive constants $A, B$ and $C$ such that
\begin{equation}\label{LW}
{\mathbb P}\Big(\big|W_1(\mu_n, \mu)- {\mathbb E}[W_1(\mu_n, \mu)]\big|>x\Big) \leq \begin{cases}
  2\exp\left(- n A x^p   \right) \quad  \text{if\ \ $x \geq C$}\\
   2\exp\left(- nB x^2  \right) \quad   \text{if\ \ $x \in [0, C]$.}
  \end{cases}
\end{equation}
This follows from Proposition \ref{cram} (case $p=1$) and  Proposition \ref{vprod} (case $p>1$).
\item If \eref{laplace2} holds for some $p \in (0,1)$, then
there exist some positive constants $A, B, C, D$ and $L$ such that
$$
{\mathbb P}\Big(\big|W_1(\mu_n, \mu)- {\mathbb E}[W_1(\mu_n, \mu)]\big|>x\Big) \leq \begin{cases}
  C\exp\left(- n^p A x^p   \right) \quad  \text{if\ \ $x \geq L n^{-(1-p)/(2-p)}$}\\
   D\exp\left(- nB x^2  \right) \quad   \text{if\ \ $x \in [0, L n^{-(1-p)/(2-p)}]$.}
  \end{cases}
$$
This follows from Proposition \ref{findsa}.
\item If \eref{weakVB} holds for some  $p \in (1,2)$, then there exists a positive
constant $C$ such that
$$
{\mathbb P}\Big(\big|W_1(\mu_n, \mu)- {\mathbb E}[W_1(\mu_n, \mu)]\big|>x\Big) \leq \frac{C}{n^{p-1} x^p} \, .
$$
This follows from Proposition \ref{weakVBEI}.
\item If \eref{weakVB} holds for some $p\geq 2$, then there exists a positive constant $C$
such that
$$
{\mathbb P}\Big(\big|W_1(\mu_n, \mu)- {\mathbb E}[W_1(\mu_n, \mu)]\big|>x\Big) \leq \frac{C}{
n^{p/2} x^p} \, .
$$
This follows from Proposition \ref{RCG}.
\end{itemize}
And for the moment bounds of $S_n$:
\begin{itemize}
\item  If \eref{fnf58} for some $p \in [1,2]$, then
\begin{equation}\label{MZ}
\Big\|W_1(\mu_n, \mu)- {\mathbb E}[W_1(\mu_n, \mu)]\Big\|^p_p \leq \frac{C}{ n^{p-1}}\, .
\end{equation}
This follows from Proposition \ref{VBEI}.
\item  If \eref{shnbs} holds  for some $p\geq 2$, then
\begin{equation}\label{ROSEN}
\Big\|W_1(\mu_n, \mu)- {\mathbb E}[W_1(\mu_n, \mu)]\Big\|^p_p \leq \frac{C}{ n^{p/2}}\, .
\end{equation}
This follows from Proposition \ref{Ros}.
\end{itemize}
Let us now give some references on the subject.

As already mentioned, the subgaussian bound \eref{LW} for $p=2$ is proved in
the paper by Djellout {\it et al.} \cite{DGW}. Notice that these authors also
consider the Wasserstein metrics $W_r$ for $r\geq 1$, with cost function
$c(x,y)=(d(x,y))^r$.

 In the iid case, when $X_i=\varepsilon_i$,
some very precise results are given in the paper by Gozlan and Leonard \cite{GL}, for a more general class
of Wasserstein metrics (meaning that the cost function is not necessary a distance). In
the case of $W_1$, they have obtained deviation inequalities under some conditions of the
Laplace transform of some convex and increasing function of $d(x_0,X_1)$ (see their Theorem 10 combined with their Theorem 7). In particular, {\it via} their Lemma 1, they have obtained a Cram\'er-type inequality
for $W_1$ similar to what we get in Proposition \ref{cram}.

In the dependent case, another important reference is the recent paper by Chazottes and Gou\"ezel \cite{CG12}.
These authors consider separately Lipschitz functionals of iterates of maps that
can be modeled by Young towers. They obtain exponential or polynomial bounds according as
the covariances between Lipschitz functions of the iterates decrease with an exponential or polynomial rate. See their Section 7.3  for the applications to the Wassertein distance $W_1$. Note that the Markov chains associated to the maps considered by Chazottes and Gou\"ezel do not in general
satisfy the one step contraction, and are much more difficult to handle than the class
of Markov chains of the present paper.

\subsection{Discussion}
Of course, the next question is that of the behavior of ${\mathbb E}[W_1(\mu_n, \mu)]$,
because it can give us information on $W_1(\mu_n, \mu)$ through the preceding inequalities.
For instance, from \eref{MZ}, we infer that if \eref{fnf58} holds for some $p\in [1,2]$, then
\begin{equation}\label{clearbis}
{\mathbb E}[W_1(\mu_n, \mu)] \leq \|W_1(\mu_n, \mu)\|_p \leq {\mathbb E}[W_1(\mu_n, \mu)] + \frac {C}{ n^{(p-1)/p} }\, .
\end{equation}
In the same way, from \eref{ROSEN}, we infer  that if \eref{shnbs} holds  for some $p\geq 2$, then
\begin{equation}\label{clear}
{\mathbb E}[W_1(\mu_n, \mu)] \leq \|W_1(\mu_n, \mu)\|_p \leq {\mathbb E}[W_1(\mu_n, \mu)] + \frac {C}{ \sqrt n} \, .
\end{equation}

Let  us first quote that, if ${\mathbb E}[G_{X_1}(X_1)] <\infty$ and ${\mathbb E}[G_\varepsilon(\varepsilon)] <\infty$, then
${\mathbb E}[W_1(\mu_n, \mu)]$ converges to $0$. Indeed, the Markov chain $(X_i)_{i \geq 1}$ satisfies the strong law of large numbers:
$$
  \lim_{n \rightarrow \infty} \mu_n(f)= \mu(f) \quad \text{almost surely,}
$$
for any $f$ such that $f(x) \leq C(1+d(x_0,x))$. Hence,
it follows from Theorem 6.9 in Villani \cite{V} that $W_1(\mu_n, \mu)$
converges to $0$ almost surely, and that ${\mathbb E}[W_1(\mu_n, \mu)]$ converges to $0$.

The question of the rate of convergence to 0 of ${\mathbb E}[W_1(\mu_n, \mu)]$ is delicate, and has a long history.
Let us recall some know results in the iid case, when $X_i=\varepsilon_i$.
\begin{itemize}
\item If ${\mathcal X}={\mathbb R}$ and $d(x,y)=|x-y|$, and if $\int |x|\sqrt{{\mathbb P}(|X_1|>x)} dx < \infty$, then
$$ \lim_{n \rightarrow \infty} \sqrt{n}{\mathbb E}[W_1(\mu_n, \mu)]=c$$
with $c\neq 0$ as soon as $X_i$ is not almost surely constant.
This follows from del Barrio {\it et al.} \cite{dBarG} and can be easily extended
to our Markov setting.

\item If ${\mathcal X}={\mathbb R}^\ell$ and $d(x,y)=\|x-y\|$ for some norm $\|\cdot\|$, let us recall some recent results by  Fournier and Guillin \cite{FG}
(see also Dereich {\it et al.} \cite{DSS}). In Theorem 1 of Fournier and Guillin \cite{FG}, the following upper bounds are proved:
Assume that $p>1$ and that $\int \|x\|^p \mu(dx) <\infty$, then
\begin{equation}\label{dereich}
{\mathbb E}[W_1(\mu_n, \mu)]\leq
\begin{cases}
  C (n^{-1/2} + n^{-(p-1)/p}) \quad \text{if $\ell=1$ and $p \neq 2$}\\
  C (n^{-1/2}\ln(1+n) + n^{-(p-1)/p}) \quad   \text{if $\ell=2$ and $p \neq 2$}\\
  C (n^{-1/\ell} + n^{-(p-1)/p}) \quad \text{if $\ell>2$ and $p \neq \ell/(\ell-1)$.}
  \end{cases}
\end{equation}
Combining this upper bound with \eref{clearbis} and \eref{clear}, we
obtain the following proposition
\begin{prop} Let $X_1, \ldots, X_n$ be an iid sequence
of ${\mathbb R}^\ell$-valued random variables, with common distribution $\mu$.
let $p>1$ and assume  that $\int \|x\|^p \mu(dx) <\infty$.
Then the quantity $\|W_1(\mu_n, \mu)\|_p$ satifies the upper bound \eref{dereich}.
\end{prop}
Note that Fournier and Guillin \cite{FG}  consider the case
of $W_r$  metrics, and the upper bound \eref{dereich}
is just a particular case of their Theorem 1.
Note also that an extension of inequality \eref{dereich}
to $\rho$-mixing Markov chains is given in Theorem 15
of the same paper.

 In their Theorem 2, Fournier and Guillin \cite{FG} give some deviation inequalities for $${\mathbb P}\big(W_r(\mu_n, \mu)>x\big).$$
 For $r=1$, these results are different from ours, since they do not deal with concentration around the mean.
In particular their upper bounds depend on the dimension $\ell$,
and  for $r=1$ and $\ell \geq 3$ they are useless for $x= yn^{-\alpha}$ as soon as
$\alpha \in (1/\ell, 1/2]$. This is coherent with our upper bounds
of Section \ref{DU} since in that case ${\mathbb E}[W_1(\mu_n, \mu)]$
can be of order $n^{-1/\ell}$. Let us note,  however,
 that the results
of Section \ref{DU}  give always an efficient upper bound for the concentration of $ W_1(\mu_n, \mu)$ around ${\mathbb E}[W_1(\mu_n, \mu)]$ for
any $x= yn^{-\alpha}$ with $\alpha \in [0, 1/2]$, that is in
the whole range from small to large deviations, whatever
the dimension of ${\mathcal X}$.

\item Concerning the behavior of ${\mathbb E}[W_1(\mu_n, \mu)]$ in the infinite dimensional case, let us mention the upper bound (15) in  Boissard
\cite{Boi}. This upper bound involves the covering numbers of an increasing sequence of
compact sets $K_t$ for which $\mu(K^c_t)$ tends to zero as $t$ tends to infinity.
Some extensions to a class of Markov chains are given in Section 2.4 of the same paper.
In particular, his results apply to one step contracting Markov chains satisfying \eref{contract} (again, this follows from Proposition 3.1 of Djellout {\it et al.}
\cite{DGW}).

\end{itemize}

\bigskip

\noindent \textbf{Acknowledgements.} J\'er\^ome Dedecker is partially supported by the French ANR project TopData.

\end{document}